\documentclass{article}

\usepackage{amsmath,amssymb,amsthm}
\usepackage{a4wide}
\usepackage[all]{xy}
\usepackage{subfig}
\usepackage{pdfsync}
\usepackage{graphicx}
\usepackage{url}
\usepackage[usenames,dvipsnames]{pstricks}
\usepackage{epsfig}
\usepackage{pst-grad}
\usepackage{pst-plot}

\newtheorem{theorem}{Theorem}[section]
\newtheorem{lemma}[theorem]{Lemma}

\def\R{\mathbb{R}}

\renewcommand{\geq}{\geqslant}
\renewcommand{\leq}{\leqslant}

\begin{document}

\title{Modeling TB-HIV syndemic and treatment\footnote{This is a preprint of a paper 
whose final and definite form is: Journal of Applied Mathematics (ISSN 1110-757X) 2014, 
Article ID 248407, {\tt http://dx.doi.org/10.1155/2014/248407}}}

\author{Cristiana J. Silva\\
\texttt{cjoaosilva@ua.pt}
\and Delfim F. M. Torres\footnote{Corresponding author.}\\
\texttt{delfim@ua.pt}}

\date{Center for Research and Development in Mathematics and Applications (CIDMA)\\
Department of Mathematics, University of Aveiro, 3810-193 Aveiro, Portugal}

\maketitle


\begin{abstract}
Tuberculosis (TB) and human immunodeficiency virus (HIV) can be considered
a deadly human syndemic. In this article, we formulate a model for TB
and HIV transmission dynamics. The model considers both TB and acquired immune
deficiency syndrome (AIDS) treatment for individuals with only one of the
infectious diseases or both. The basic reproduction number and equilibrium
points are determined and stability is analyzed. Through simulations, we show
that TB treatment for individuals with only TB infection reduces the number
of individuals that become co-infected with TB and HIV/AIDS, and reduces
the diseases (TB and AIDS) induced deaths. Analogously, the treatment
of individuals with only AIDS also reduces the number of co-infected individuals.
Further, TB-treatment for co-infected individuals in the active and latent
stage of TB disease, implies a decrease of the number of individuals
that passes from HIV-positive to AIDS.
\end{abstract}


\paragraph{Keywords:} Tuberculosis, Human immunodeficiency virus,
Syndemic, Treatment, Equilibrium, Stability.

\paragraph{Mathematics Subject Classification 2010:} 34D30; 92D30; 93A30.


\section{Introduction}

Tuberculosis (TB) and human immunodeficiency virus/acquired immune deficiency syndrome (HIV/AIDS)
are the leading causes of death from an infectious disease worldwide \cite{TB_WHO_report_2013}.
Individuals infected with HIV are more likely to develop TB disease because of their immunodeficiency,
and HIV infection is the most powerful risk factor for progression from TB infection to disease
\cite{Getahun:etall:CID2010}. This interaction justifies the fact that HIV and TB can be considered
a deadly human \emph{syndemic}, where syndemic refers to the convergence of two or more diseases
that act synergistically to magnify the burden of disease \cite{Kwan_Ernst_HIV_TB_Syndemic}.

Following UNAIDS global report on AIDS epidemic 2013 \cite{UNAIDS_report_2013}, globally,
an estimated 35.3 million people were living with HIV in 2012. An increase from previous years,
as more people are receiving the life-saving antiretroviral therapy (ART). There were approximately
2.3 million new HIV infections globally, showing a 33\% decline in the number of new infections
with respect to 2001. At the same time, the number of AIDS deaths is also declining with around 1.6
million AIDS deaths in 2012, down from about 2.3 million in 2005.
In 2012, 1.1 million of 8.6 million people who developed TB worldwide were HIV-positive. The number
of people dying from HIV-associated TB has been falling since 2003. However, there were still
320 000 deaths from HIV-associated TB in 2012 and further efforts are needed to reduce this burden
\cite{TB_WHO_report_2013}. ART is a critical intervention for reducing the risk of TB morbidity
and mortality among people living with HIV and, when combined with isoniazid preventive therapy,
it can have a significant impact on TB prevention \cite{TB_WHO_report_2013}.

Collaborative TB/HIV activities (including HIV testing, ART therapy and TB preventive measures)
are crucial for the reduction of TB-HIV coinfected individuals. The World Health Organization (WHO)
estimates that these collaborative activities prevented 1.3 million people from dying,
from 2005 to 2012. However, significant challenges remain: the reduction of tuberculosis
related deaths among people living with HIV has slowed in recent years; the ART therapy
is not being delivered to TB-HIV coinfected patients in the majority of the countries
with the largest number of TB/HIV patients; the pace of treatment scale-up for TB/HIV
patients has slowed; less than half of notified TB patients were tested for HIV in 2012;
and only a small fraction of TB/HIV infected individuals received TB preventive therapy
\cite{UNAIDS_report_2013}.

The study of the joint dynamics of TB and HIV present formidable mathematical challenges
due to the fact that the models of transmission are quite distinct \cite{CChavez_TB_HIV_2009}.
Few mathematical models have been proposed for TB-HIV coinfection (see, for example,
\cite{Bhunu:BMB:2009:HIV:TB,Kirschner:TB:HIV:1999,Naresh:TB:HIV:2005,CChavez_TB_HIV_2009,Song:TB:HIV:2008}).
Kirschner \cite{Kirschner:TB:HIV:1999} developed a cellular model for HIV-1 and TB coinfection inside a host.
Roeger et al. \cite{CChavez_TB_HIV_2009} proposed a population model for TB-HIV/AIDS coinfection
transmission dynamics, assuming that TB-infected individuals in the active stage of the disease
are too ill to remain sexually active and therefore they are unable to transmit HIV. In this work
we assume that active TB-infected individuals are susceptible to HIV-infection.
Naresh and Tripathi \cite{Naresh:TB:HIV:2005} proposed a model for TB-HIV coinfection
in a variable size population with only TB treatment. Here we consider TB and HIV treatment
in different stages of the disease. Bhunu et al. \cite{Bhunu:BMB:2009:HIV:TB}
studied a TB-HIV coinfection model with both TB and HIV treatment. The authors did not
take into account that an individual co-infected with TB and HIV can effectively recover
from TB infection. We assume that TB can be cured, even in HIV-positive individuals
\cite{TB_WHO_report_2013}. Sharomi et al. \cite{Song:TB:HIV:2008} also considered
these assumptions, subdividing the total population into 15 classes. It is our aim
in this work to develop a model that balances two goals: simplicity and useful information.

The paper is organized as follows. Section~\ref{sec:2} describes our model
for TB-HIV syndemic with TB and HIV treatment. In Section~\ref{sec:3} the
positivity and boundedness of solutions of the model are proved
and in Section~\ref{sec:4} equilibrium points and respective stability are analyzed.
Section~\ref{sect:NumSim} is devoted to numerical simulations and discussion of results.


\section{TB-HIV/AIDS model}
\label{sec:2}

The model subdivides the human population into 10 mutually-exclusive compartments,
namely susceptible individuals ($S$), TB-latently infected individuals, who have
no symptoms of TB disease and are not infectious ($L_T$), TB-infected individuals,
who have active TB disease and are infectious ($I_T$), TB-recovered individuals
($R_T$), HIV-infected individuals with no clinical symptoms of AIDS ($I_H$),
HIV-infected individuals with AIDS clinical symptoms ($A$), TB-latent individuals
co-infected with HIV (pre-AIDS) ($L_{TH}$), HIV-infected individuals (pre-AIDS)
co-infected with active TB disease ($I_{TH}$), TB-recovered individuals with
HIV-infection without AIDS symptoms ($R_{TH}$), HIV-infected individuals
with AIDS symptoms co-infected with TB ($A_T$).
The total population at time $t$, denoted by $N(t)$, is given by
\begin{equation*}
N(t) = S(t) + L_T(t) + I_T(t) + R_T(t) + I_H(t) + A(t)
+ I_{TH}(t) + L_{TH}(t) + R_{TH}(t) + A_T(t)  \, .
\end{equation*}
The susceptible population is increased by the recruitment of individuals
(assumed susceptible) into the population, at a rate $\Lambda$.
All individuals suffer from natural death, at a constant rate $\mu$.
Susceptible individuals acquire TB infection from individuals
with active TB at a rate $\lambda_T$, given by
\begin{equation*}
\lambda_T = \frac{\beta_1}{N} \left(I_T + I_{TH} + A_T\right) \, ,
\end{equation*}
where $\beta_1$ is the effective contact rate for TB infection.
Similarly, susceptible individuals acquire HIV infection, following
effective contact with people infected with HIV at a rate $\lambda_H$, given by
\begin{equation*}
\lambda_H = \frac{\beta_2}{N} \left[ I_H + I_{TH} + L_{TH}
+ R_{TH} + \eta \left(A + A_T\right) \right] \, ,
\end{equation*}
where $\beta_2$ is the effective contact rate for HIV transmission
and the modification parameter $\eta \geq 1$ accounts for the relative
infectiousness of individuals with AIDS symptoms, in comparison to those
infected with HIV with no AIDS symptoms. Individuals with AIDS symptoms
are more infectious than HIV-infected individuals (pre-AIDS) because
they have a higher viral load and there is a positive correlation
between viral load and infectiousness \cite{art:viral:load}.

Individuals leave the latent-TB class $L_T$ by becoming infectious,
at a rate $k_1$, or recovered, with a treatment rate $\tau_1$.
The treatment rate for active TB-infected individuals is $\tau_2$.
We assume that TB-recovered individuals $R_T$ acquire partial immunity
and the transmission rate for this class is given by $\beta'_1 \lambda_T$
with $\beta'_1 \leq 1$. Individuals with active TB disease suffer
induced death at a rate $d_T$. We assume that individuals in the class
$R_T$ are susceptible to HIV infection at a rate $\lambda_H$.
On the other hand, TB-active infected individuals $I_T$ are susceptible
to HIV infection, at a rate $\delta \lambda_H$, where the modification
parameter $\delta \geq 1$ accounts for higher probability of individuals
in class $I_T$ to become HIV-positive.

HIV-infected individuals (with no AIDS symptoms) progress to the AIDS class $A$,
at a rate $\rho_1$. HIV-infected individuals with AIDS symptoms are treated
for HIV at the rate $\alpha_1$ and suffer induced death at a rate $d_A$.
Individuals in the class $I_H$ are susceptible to TB infection at a rate
$\psi \lambda_T$, where $\psi \geq 1$ is a modification parameter
traducing the fact that HIV infection is a driver of TB epidemic
\cite{Kwan_Ernst_HIV_TB_Syndemic}.

HIV-infected individuals (pre-AIDS) co-infected with TB-disease,
in the active stage $I_{TH}$, are treated for TB at the rate $\tau_3$
and progress to the AIDS-TB co-infection class $A_T$ at a rate $\rho_2$.
Individuals in the class $I_{TH}$ suffer TB induced death at a rate $d_T$.
The anti-TB drugs can prevent or decrease the likelihood of TB infection progression
to active TB disease in individuals in the class $L_{TH}$ \cite{USAID:TB:HIV}.
The treatment rate for individuals in this class is given by $\tau_4$.
However, individuals in the class $L_{TH}$ are more likely to progress
to active TB disease than individuals infected only with latent TB.
In our model, this progression rate is given by $k_2$.
Similarly, HIV infection makes individuals more susceptible to TB reinfection
when compared with non HIV-positive patients. The modification parameter associated
to the TB reinfection rate, for individuals in the class $R_{TH}$, is given by $\beta'_2$,
where $\beta'_2 \geq 1$. Individuals in this class progress to class $A_T$, at a rate $\rho_3$.

HIV-infected individuals (with AIDS symptoms), co-infected with TB, are treated for HIV,
at a rate $\alpha_2$. Individuals in the class $A_T$ suffer from AIDS-TB
coinfection induced death rate, at a rate $d_{TA}$.

\begin{figure}
\centering
\scalebox{0.90}
{
\begin{pspicture}(0,-5.098906)(16.862812,5.078906)
\psframe[linewidth=0.04,dimen=outer](9.280937,4.538906)(7.7009373,3.7589064)
\psframe[linewidth=0.04,dimen=outer](11.320937,-0.86109376)(9.740937,-1.6410937)
\psframe[linewidth=0.04,linestyle=dashed,dash=0.16cm 0.16cm,dimen=outer](7.2809377,0.13890626)(5.7009373,-0.64109373)
\psframe[linewidth=0.068,linestyle=dotted,dotsep=0.16cm,dimen=outer](15.300938,0.13890626)(13.720938,-0.64109373)
\psframe[linewidth=0.068,linestyle=dotted,dotsep=0.16cm,dimen=outer](13.300938,2.1389062)(11.720938,1.3589063)
\psframe[linewidth=0.04,linestyle=dashed,dash=0.16cm 0.16cm,dimen=outer](3.1209376,0.13890626)(1.5409375,-0.64109373)
\psframe[linewidth=0.04,linestyle=dashed,dash=0.16cm 0.16cm,dimen=outer](5.3609376,2.1389062)(3.7809374,1.3589063)
\psframe[linewidth=0.04,dimen=outer](7.5209374,-2.4610937)(5.9409375,-3.2410936)
\psframe[linewidth=0.04,dimen=outer](13.700937,-2.4610937)(12.120937,-3.2410936)
\psframe[linewidth=0.04,dimen=outer](10.100938,-4.1210938)(8.520938,-4.901094)
\psline[linewidth=0.04cm,arrowsize=0.05291667cm 2.0,arrowlength=1.4,arrowinset=0.4]{->}(7.6609373,4.1989064)(4.7409377,2.1989062)
\psline[linewidth=0.04cm,arrowsize=0.05291667cm 2.0,arrowlength=1.4,arrowinset=0.4]{->}(9.320937,4.1989064)(12.600938,2.2189062)
\psline[linewidth=0.04cm,arrowsize=0.05291667cm 2.0,arrowlength=1.4,arrowinset=0.4]{<-}(5.6209373,-0.28109375)(3.2409375,-0.28109375)
\psline[linewidth=0.04cm,arrowsize=0.05291667cm 2.0,arrowlength=1.4,arrowinset=0.4]{->}(7.3609376,-0.18109375)(11.660937,1.7389063)
\psline[linewidth=0.04cm,arrowsize=0.05291667cm 2.0,arrowlength=1.4,arrowinset=0.4]{->}(11.700937,1.2989062)(6.8009377,-2.4010937)
\psline[linewidth=0.04cm,arrowsize=0.05291667cm 2.0,arrowlength=1.4,arrowinset=0.4]{->}(14.260938,0.19890624)(13.300938,1.2589062)
\psline[linewidth=0.04cm,arrowsize=0.05291667cm 2.0,arrowlength=1.4,arrowinset=0.4]{->}(12.820937,1.2989062)(13.880938,0.19890624)
\psline[linewidth=0.04cm,arrowsize=0.05291667cm 2.0,arrowlength=1.4,arrowinset=0.4]{->}(9.680938,-1.2810937)(7.3209376,-2.4210937)
\psline[linewidth=0.04cm,arrowsize=0.05291667cm 2.0,arrowlength=1.4,arrowinset=0.4]{->}(13.000937,-2.3810937)(11.400937,-1.3010937)
\psline[linewidth=0.04cm,arrowsize=0.05291667cm 2.0,arrowlength=1.4,arrowinset=0.4]{->}(7.6009374,-2.8610938)(12.040937,-2.8610938)
\psline[linewidth=0.04cm,arrowsize=0.05291667cm 2.0,arrowlength=1.4,arrowinset=0.4]{->}(12.9609375,-3.2810938)(10.160937,-4.461094)
\psline[linewidth=0.04cm,arrowsize=0.05291667cm 2.0,arrowlength=1.4,arrowinset=0.4]{->}(8.480938,-4.2010937)(6.8609376,-3.3210938)
\psline[linewidth=0.04cm,arrowsize=0.05291667cm 2.0,arrowlength=1.4,arrowinset=0.4]{->}(6.2009373,-3.3410938)(8.420938,-4.601094)
\usefont{T1}{ptm}{m}{n}
\rput(8.518281,4.1539063){\large $S$}
\usefont{T1}{ptm}{m}{n}
\rput(4.668281,1.7539062){\large $L_T$}
\usefont{T1}{ptm}{m}{n}
\rput(2.3282812,-0.28609374){\large $I_T$}
\usefont{T1}{ptm}{m}{n}
\rput(12.558281,1.7339063){\large $I_H$}
\usefont{T1}{ptm}{m}{n}
\rput(6.508281,-0.24609375){\large $R_T$}
\usefont{T1}{ptm}{m}{n}
\rput(14.4582815,-0.28609374){\large $A$}
\usefont{T1}{ptm}{m}{n}
\rput(10.638281,-1.2460938){\large $L_{TH}$}
\usefont{T1}{ptm}{m}{n}
\rput(6.918281,-2.8660936){\large $I_{TH}$}
\usefont{T1}{ptm}{m}{n}
\rput(13.058281,-2.8660936){\large $R_{TH}$}
\usefont{T1}{ptm}{m}{n}
\rput(9.308281,-4.526094){\large $A_{T}$}
\usefont{T1}{ptm}{m}{n}
\rput{34.0}(2.9529414,-2.8401985){\rput(6.1023436,3.4289062){$\lambda_T$}}
\usefont{T1}{ptm}{m}{n}
\rput{-30.0}(-0.21864277,6.0024514){\rput(11.072344,3.4289062){$\lambda_H$}}
\usefont{T1}{ptm}{m}{n}
\rput{-45.0}(3.511642,10.136291){\rput(13.972343,0.8489063){$\alpha_1$}}
\usefont{T1}{ptm}{m}{n}
\rput{-50.0}(4.260766,10.31568){\rput(13.172344,0.60890627){$\rho_1$}}
\usefont{T1}{ptm}{m}{n}
\rput{-28.0}(2.6423028,3.3361356){\rput(7.992344,-3.6110938){$\alpha_2$}}
\usefont{T1}{ptm}{m}{n}
\rput{-30.0}(3.0368948,3.0722833){\rput(7.2323437,-4.1110935){$\rho_2$}}
\usefont{T1}{ptm}{m}{n}
\rput{23.0}(-0.6896543,-4.7918034){\rput(11.412344,-4.0710936){$\rho_3$}}
\usefont{T1}{ptm}{m}{n}
\rput{23.0}(0.9393141,-3.638439){\rput(9.3923435,0.50890625){$\lambda_H$}}
\usefont{T1}{ptm}{m}{n}
\rput{37.0}(1.5412767,-5.507953){\rput(8.982344,-0.43109375){$\psi \lambda_T$}}
\usefont{T1}{ptm}{m}{n}
\rput(9.852344,-3.0510938){$\tau_3$}
\usefont{T1}{ptm}{m}{n}
\rput{25.0}(0.141107,-3.778055){\rput(8.572344,-1.5510937){$k_2$}}
\usefont{T1}{ptm}{m}{n}
\rput{-35.0}(3.1449733,6.753019){\rput(12.262343,-1.5910938){$\beta^{'}_2 \lambda_T$}}
\usefont{T1}{ptm}{m}{n}
\rput(8.782344,4.808906){$\Lambda$}
\psline[linewidth=0.04cm,arrowsize=0.05291667cm 2.0,arrowlength=1.4,arrowinset=0.4]{->}(9.340938,4.3189063)(10.120937,4.3389063)
\usefont{T1}{ptm}{m}{n}
\rput(9.582344,4.4889064){$\mu$}
\usefont{T1}{ptm}{m}{n}
\rput(13.6423435,1.9489063){$\mu$}
\usefont{T1}{ptm}{m}{n}
\rput(14.1423435,-3.0510938){$\mu$}
\usefont{T1}{ptm}{m}{n}
\rput(7.722344,-0.51109374){$\mu$}
\psline[linewidth=0.04cm,arrowsize=0.05291667cm 2.0,arrowlength=1.4,arrowinset=0.4]{->}(5.9009376,-3.0810938)(5.1809373,-3.0810938)
\usefont{T1}{ptm}{m}{n}
\rput(0.82234377,-0.49109375){$\mu+d_T$}
\usefont{T1}{ptm}{m}{n}
\rput(5.4823437,-3.3310938){$\mu+d_T$}
\usefont{T1}{ptm}{m}{n}
\rput(10.882343,-4.8710938){$\mu+d_{TA}$}
\usefont{T1}{ptm}{m}{n}
\rput(3.2423437,1.5689063){$\mu$}
\usefont{T1}{ptm}{m}{n}
\rput(11.822344,-0.95109373){$\mu$}
\usefont{T1}{ptm}{m}{n}
\rput(15.962344,-0.43109375){$\mu+d_A$}
\usefont{T1}{ptm}{m}{n}
\rput(4.3123436,-0.49109375){$\tau_2$}
\psline[linewidth=0.04cm,arrowsize=0.05291667cm 2.0,arrowlength=1.4,arrowinset=0.4]{->}(1.5009375,-0.26109374)(0.5209375,-0.26109374)
\psline[linewidth=0.04cm,arrowsize=0.05291667cm 2.0,arrowlength=1.4,arrowinset=0.4]{->}(7.3409376,-0.26109374)(8.220938,-0.26109374)
\psline[linewidth=0.04cm,arrowsize=0.05291667cm 2.0,arrowlength=1.4,arrowinset=0.4]{->}(10.120937,-4.601094)(11.200937,-4.601094)
\psline[linewidth=0.04cm,arrowsize=0.05291667cm 2.0,arrowlength=1.4,arrowinset=0.4]{->}(2.3009374,-0.70109373)(5.8809376,-2.8010938)
\psline[linewidth=0.04cm,arrowsize=0.05291667cm 2.0,arrowlength=1.4,arrowinset=0.4]{->}(4.3009377,1.3389063)(2.3409376,0.19890624)
\psline[linewidth=0.04cm,arrowsize=0.05291667cm 2.0,arrowlength=1.4,arrowinset=0.4]{<-}(4.9409375,1.3189063)(6.5009375,0.23890625)
\usefont{T1}{ptm}{m}{n}
\rput{-34.0}(0.60520494,3.5579736){\rput(6.1023436,0.80890626){$\beta^{'}_1 \lambda_T$}}
\usefont{T1}{ptm}{m}{n}
\rput{32.0}(0.9818715,-1.5257553){\rput(3.1323438,0.9689062){$k_1$}}
\psline[linewidth=0.04cm,arrowsize=0.05291667cm 2.0,arrowlength=1.4,arrowinset=0.4]{->}(4.7209377,1.2789062)(6.2809377,0.19890624)
\usefont{T1}{ptm}{m}{n}
\rput{-34.0}(0.6034248,3.1121511){\rput(5.3723435,0.5889062){$\tau_1$}}
\usefont{T1}{ptm}{m}{n}
\rput{-30.0}(1.4413812,1.6777455){\rput(3.8323438,-1.8310938){$\delta \lambda_H$}}
\psline[linewidth=0.04cm,arrowsize=0.05291667cm 2.0,arrowlength=1.4,arrowinset=0.4]{->}(11.420938,-1.1210938)(12.480938,-1.1210938)
\psline[linewidth=0.04cm,arrowsize=0.05291667cm 2.0,arrowlength=1.4,arrowinset=0.4]{->}(15.340938,-0.18109375)(16.220938,-0.18109375)
\psline[linewidth=0.04cm,arrowsize=0.05291667cm 2.0,arrowlength=1.4,arrowinset=0.4]{->}(11.380938,-1.6010938)(12.560938,-2.3810937)
\usefont{T1}{ptm}{m}{n}
\rput{-32.0}(2.8905547,5.899){\rput(11.712344,-2.0710938){$\tau_4$}}
\psline[linewidth=0.04cm,arrowsize=0.05291667cm 2.0,arrowlength=1.4,arrowinset=0.4]{->}(13.760938,-2.8810937)(14.540937,-2.8810937)
\psline[linewidth=0.04cm,arrowsize=0.05291667cm 2.0,arrowlength=1.4,arrowinset=0.4]{->}(13.380938,1.7389063)(14.160937,1.7389063)
\psline[linewidth=0.04cm,arrowsize=0.05291667cm 2.0,arrowlength=1.4,arrowinset=0.4]{<-}(2.9209375,1.7389063)(3.7009375,1.7389063)
\psline[linewidth=0.04cm,arrowsize=0.05291667cm 2.0,arrowlength=1.4,arrowinset=0.4]{->}(8.520938,5.058906)(8.520938,4.538906)
\end{pspicture}}
\caption{Model for TB-HIV/AIDS transmission with treatment.}
\label{fig:model:flow}
\end{figure}
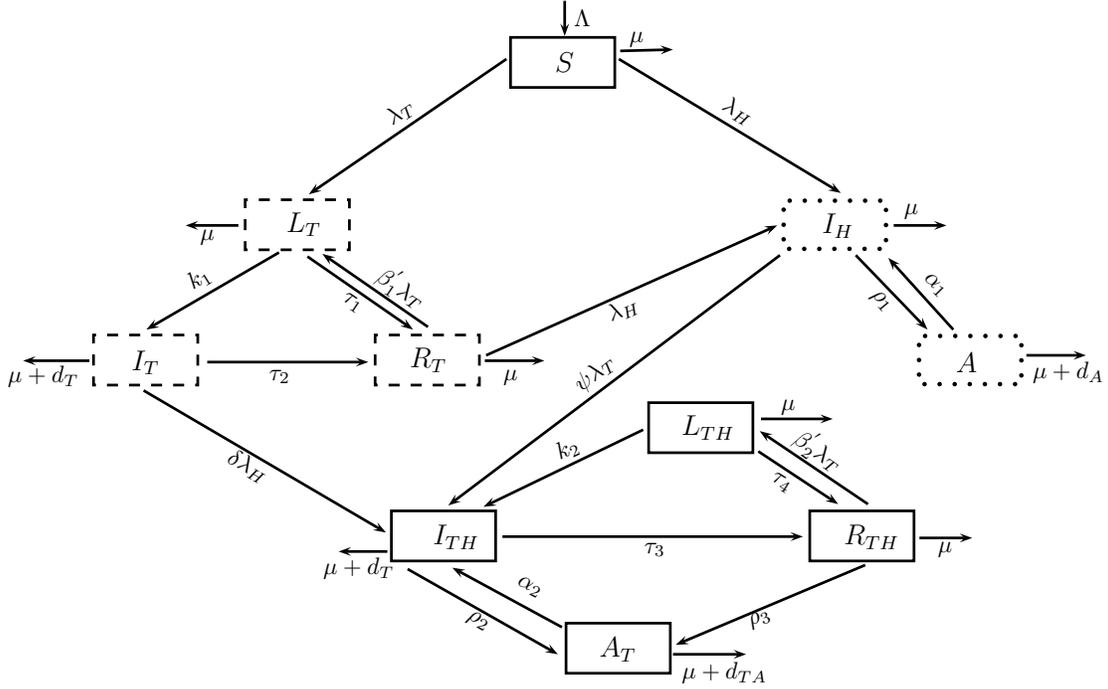

The aforementioned assumptions result in the following system of differential
equations that describes the transmission dynamics of TB and HIV disease:
\begin{equation}
\label{model:TB:HIV}
\begin{cases}
\dot{S}(t) = \Lambda - \lambda_T S(t) - \lambda_H S(t) - \mu S(t),\\[0.2 cm]
\dot{L}_T(t) = \lambda_T S(t) + \beta^{'}_1 \lambda_T R_T(t) - (k_1 + \tau_1 + \mu)L_T(t),\\[0.2 cm]
\dot{I}_T(t) = k_1 L_T(t) - (\tau_2 +d_T +\mu + \delta \lambda_H)I_T(t), \\[0.2 cm]
\dot{R}_T(t) = \tau_1 L_T(t) + \tau_2 I_T(t) - (\beta^{'}_1 \lambda_T + \lambda_H + \mu) R_T(t),\\[0.2 cm]
\dot{I}_H(t) = \lambda_H S(t) - (\rho_1 + \psi \lambda_T + \mu)I_H(t) + \alpha_1 A(t) + \lambda_H R_T(t), \\[0.2 cm]
\dot{A}(t) =  \rho_1 I_H(t) - \alpha_1 A(t) - (\mu + d_A) A(t),\\[0.2 cm]
\dot{L}_{TH}(t) = \beta^{'}_2 \lambda_T R_{TH}(t) - (k_2 + \tau_4 + \mu) L_{TH}(t),\\[0.2 cm]
\dot{I}_{TH}(t) = \delta \lambda_H I_T(t) + \psi \lambda_T I_H(t) + \alpha_2 A_T(t)+ k_2 L_{TH}(t)
- (\tau_3 + \rho_2 + \mu + d_T)I_{TH}(t),\\[0.2 cm]
\dot{R}_{TH}(t) = \tau_3 I_{TH}(t) + \tau_4 L_{TH}(t) - (\beta^{'}_2 \lambda_T  + \rho_3 + \mu)R_{TH}, \\[0.2 cm]
\dot{A}_T(t) = \rho_2 I_{TH}(t) + \rho_3 R_{TH} -(\alpha_2 + \mu + d_{TA})A_T(t) \, .
\end{cases}
\end{equation}
The model flow is described in Figure~\ref{fig:model:flow}.
The initial conditions of model \eqref{model:TB:HIV} satisfy
\begin{equation}
\label{eq:init:cond:geral}
\begin{split}
&S(0)= S_0 \geq 0 \, , \quad \quad L_T(0) = L_{T0} \geq 0 \, ,
\quad \quad I_T(0) = I_{T0} \geq 0 \, ,
\quad \quad  R_T(0) = R_{T0} \geq 0 \, ,\\
&I_H(0) = I_{H0} \geq 0 \, , \quad \quad A(0) = A_0 \geq 0 \, ,
\quad \quad L_{TH}(0) = L_{TH0} \geq 0 \, ,\\
&I_{TH}(0) = I_{TH0} \geq 0 \, , \quad \quad R_{TH}(0) = R_{TH0} \geq 0
\, , \quad \quad A_T(0) = A_{T0} \geq 0  \, .
\end{split}
\end{equation}
Note that if we consider the sub-model of \eqref{model:TB:HIV}
with no HIV/AIDS disease, that is, $I_H = A = L_{TH} = I_{TH}
= R_{TH} = A_T = 0$, then we obtain the TB model from \cite{Castillo_Chavez_1997}.
On the other hand, if we consider the sub-model with no TB, that is,
$L_T = I_T = R_T = L_{TH} = I_{TH} = R_{TH} = A_T = 0$, then we obtain an
HIV/AIDS model based on the models proposed in
\cite{Bhunu:BMB:2009:HIV:TB,Hyman:MathBio:1999}.


\section{Positivity and boundedness of solutions}
\label{sec:3}

Let $(S, L_T, I_T, R_T, I_H, A, L_{TH}, I_{TH}, R_{TH}, A_T) \in \R^{10}_{+}$
be any solution of \eqref{model:TB:HIV} with initial conditions
\eqref{eq:init:cond:geral}. Consider the biologically feasible region given by
\begin{equation}
\label{eq:feasible:region}
\Omega = \left\{ \left(S, L_T, I_T, R_T, I_H, A, L_{TH}, I_{TH}, R_{TH}, A_T\right) \,
\in \R^{10}_{+} \, : \, 0 \leq N(t) \leq \frac{\Lambda}{\mu} \right\} \, .
\end{equation}
For the model system \eqref{model:TB:HIV} to be epidemiologically meaningful,
it is important to prove that all its state variables are nonnegative
for all time $t > 0$. Suppose, for example, that at some $\bar{t} > 0$
the variable $L_T$ becomes zero, i.e., $L_T(\bar{t}) = 0$, while all other
variables are positive. Then, from the $L_T$ equation we have $d L_T(\bar{t})/dt > 0$.
Thus, $L_T(t) \geq 0$ for all $t > 0$. Analogously, we can prove
that all variables remain nonnegative for all time $t > 0$.

Adding all equations in model \eqref{model:TB:HIV} gives
\begin{equation*}
\frac{d N}{dt}(t) = \Lambda - \mu N(t) - d_T I_T(t)
- d_A A(t) - d_T I_{TH}(t) - d_{TA} A_T(t) \, .
\end{equation*}
Since $N(t) \geq I_T(t) + A(t) + I_{TH}(t) + A_T(t)$, then
\begin{equation*}
\Lambda - (\mu + d_T + d_A + d_{TA})N(t) \leq \frac{d N}{dt}(t)
\leq \Lambda - \mu N(t) \, .
\end{equation*}
Therefore, we conclude that $N(t)$ is bounded for all $t > 0$
and every solution of system \eqref{model:TB:HIV} with initial
condition in $\Omega$ remains in $\Omega$. This result is summarized below.

\begin{lemma}
The region $\Omega$ is positively invariant for the model
\eqref{model:TB:HIV} with non-negative initial conditions
in $\R^{10}_{+}$.
\end{lemma}


\section{Stability analysis}
\label{sec:4}

The model \eqref{model:TB:HIV} has four non-negative equilibria, namely
\begin{itemize}
\item[(i)] The disease-free equilibrium (no disease)
\begin{equation}
\label{eq:DFE:model:TBHIV}
\Sigma_0 = \left(S_0, L_{T_0}, I_{T_0}, R_{T_0}, I_{H_0}, A_0, L_{TH_0},
I_{TH_0}, R_{TH_0}, A_{T_0}\right) = \left(\frac{\Lambda}{\mu},
0, 0, 0, 0, 0, 0, 0, 0, 0 \right) \, .
\end{equation}

\item[(ii)] The HIV-AIDS free equilibrium
\begin{equation*}
\Sigma_T = \left(S^\diamond, L_T^\diamond, I_{T}^\diamond, R_{T}^\diamond,
I_{H}^\diamond, A^\diamond, L_{TH}^\diamond, I_{TH}^\diamond,
R_{TH}^\diamond, A_{T}^\diamond\right)
\end{equation*}
with $I_{T}^\diamond > 0$ and $I_{H}^\diamond = A^\diamond = L_{TH}^\diamond
= I_{TH}^\diamond = R_{TH}^\diamond = A_{T}^\diamond = 0$ for
$R_1 > 1$, where $R_1$ is the basic reproduction number of the model
\eqref{model:TB:HIV} with $I_H = A = L_{TH} = I_{TH} = R_{TH} = A_T = 0$
(only TB model) that is given by
\begin{equation}
\label{eq:R1}
R_1 = \frac{\Lambda}{N\mu}\left(\frac{\beta_1}{d_T+\mu+\tau_2} \right)
\left(\frac{k_1}{k_1+\tau_1+\mu } \right)
\end{equation}
(see \cite{Castillo_Chavez_1997}).

\item[(iii)] The TB-free equilibrium
\begin{equation*}
\Sigma_H = \left(S^\star, L_T^\star, I_{T}^\star, R_{T}^\star, I_{H}^\star,
A^\star, L_{TH}^\star, I_{TH}^\star, R_{TH}^\star, A_{T}^\star\right)
\end{equation*}
with
$L_T^\star = I_{T}^\star = R_{T}^\star = L_{TH}^\star
= I_{TH}^\star = R_{TH}^\star = A_{T}^\star = 0$ and
\begin{equation*}
S^\star = \frac{\Lambda}{\mu R_2}\, , \quad I_H^\star
= (R_2 - 1)\frac{\mu N_H (\alpha_1 + d_A + \mu)}{\beta_2 (\alpha_1
+ d_A + \mu + \eta \rho_1)}\, , \, \, \quad  A^\star
= (R_2 - 1) \frac{\rho_1 \mu N_H}{\beta_2 (\alpha_1 + d_A + \mu + \eta \rho_1)},
\end{equation*}
for $R_2 > 1$, where $R_2$ is the basic reproduction number of model \eqref{model:TB:HIV}
with $L_T = I_{T} = R_{T} = L_{TH} = I_{TH} = R_{TH} = A_{T} = 0$ (only HIV-AIDS model), that is,
\begin{equation}
\label{eq:R2}
R_2 =\frac{\Lambda}{N\mu} \beta_2 \left(\frac{\mu+\alpha_1+d_A
+\eta\,\rho_1 }{\mu\,\alpha_1+ (\mu + \rho_1)(\mu + d_A)}\right)\, .
\end{equation}

\item[(iv)] The syndemic equilibrium
\begin{equation*}
\Sigma^* = (S^*, L_T^*, I_T^*, R_T^*, I_H^*, A^*, L^*_{TH}, I_{TH}^*, R_{TH}^*, A_T^*)
\end{equation*}
with $I_T^* > 0$, $I_H^* > 0$, $A^* > 0$, $L_{TH}^* > 0$, $I_{TH}^* > 0$, $R_{TH}^*>0$ and  $A_T^* > 0$,
for $R_0 > 1$, where $R_0$ is the basic reproduction number of the model \eqref{model:TB:HIV}, that is,
\begin{equation*}
R_0 = \max \{ R_1, R_2 \} \, .
\end{equation*}
\end{itemize}
The details of the computation of the basic reproduction number $R_0$ are given in Appendix~\ref{A.1}.

The following theorem states the stability of the equilibrium points.

\begin{theorem}
\label{theo:stab:equil}
The disease free equilibrium $\Sigma_0$ is locally asymptotically stable
if $R_0 < 1$, and unstable if either $R_i > 1$ with $i=1, 2$.
The HIV-AIDS free equilibrium $\Sigma_T$ is locally asymptotically stable
if $R_1 > 1$, and the TB-free equilibrium $\Sigma_H$ is locally
asymptotically stable for $R_2$ near 1.
\end{theorem}

Details of the proof of Theorem~\ref{theo:stab:equil} are given in Appendix~\ref{A.2}.

\medskip

Explicit expressions for the coinfection endemic equilibrium $\Sigma^*$
are very difficult to compute analytically. In Section~\ref{sect:NumSim},
we consider an example, with $R_0 > 1$, for which there exists a syndemic
equilibrium, and analyze, numerically, the local asymptotical stability
of the syndemic equilibrium $\Sigma^*$.


\section{Numerical analysis and discussion}
\label{sect:NumSim}

For numerical simulations, we consider the following initial
conditions for system \eqref{model:TB:HIV}:
\begin{multline}
\label{init:cond}
\left(S(0), L_T(0), I_T(0), R_T(0), I_H(0), A(0), L_{TH}(0), I_{TH}(0), R_{TH}(0), A_T(0)\right)\\
= \left(\frac{60N}{100}, \frac{14 N}{100}, \frac{3N}{100}, 0, \frac{4N}{100},
\frac{N}{100}, \frac{12N}{100}, \frac{5N}{100}, 0, \frac{N}{100}\right)
\end{multline}
with $N = 50000$. The parameters of model \eqref{model:TB:HIV} take
the values of Table~\ref{table:parameters:TB-HIV}.
\begin{table}[!htb]
\centering
\begin{tabular}{l  l  l  l  l  l}
\hline
\hline
{\small{Symbol}}  &  {\small{Value}} & {\small{References}}
&  {\small{Symbol}}  &  {\small{Value}} & {\small{References}}\\
\hline
{\small{$\Lambda$}}  & {\small{$714$}} &   &
{\small{$\tau_4$}} & {\small{$1 \, yr^{-1}$}} & \\
{\small{$\mu$}} &  {\small{$1/70 \, yr^{-1}$}} &   &
{\small{$\rho_1$}} & {\small{$0.1 \, yr^{-1}$}}
& {\small{\cite{wiki:HIV:progression,site:HIV:progression:rate}}}\\
{\small{$\beta_1$}} & {\small{variable}}  &
& {\small{$\rho_2$}}  & {\small{$0.25 \, yr^{-1}$}} & \\
{\small{$\beta_2$}} & {\small{variable}} &   &
{\small{$\rho_3$}} & {\small{$0.125 \, yr^{-1}$}} & \\
{\small{$\beta'_1$}} & {\small{$0.9$}} &  &  {\small{$\alpha_1$}}
& {\small{$0.33 \, yr^{-1}$}} & {\small{\cite{Bhunu:BMB:2009:HIV:TB}}}\\
{\small{$\beta'_2$}} & {\small{$1.1$}} &   &
{\small{$\alpha_2$}} & {\small{$0.33 \, yr^{-1}$}} & \\
{\small{$k_1$}} & {\small{$1$}} & {\small{\cite{Castillo_Chavez_1997}}}
& {\small{$\psi$}} & {\small{$1.07$}} & \\
{\small{$k_2$}} & {\small{$1.3 \, k_1$}}  & \cite{USAID:TB:HIV}
& {\small{$d_T$}}  & {\small{$1/8 \, yr^{-1}$}} & \\
{\small{$\tau_1$}} & {\small{$1 \, yr^{-1}$}} &
{\small{\cite{Castillo_Chavez_1997}}}
& {\small{$d_A$}} & {\small{$0.3 \, yr^{-1}$}} & \\
{\small{$\tau_2$}} & {\small{$2 \, yr^{-1}$}}  &
{\small{\cite{Castillo_Chavez_1997}}} &
{\small{$d_{TA}$}} & {\small{$0.33 \, yr^{-1}$}} &  \\
{\small{$\tau_3$}} & {\small{$2 \, yr^{-1}$}} &  &
{\small{$\eta$}} & {\small{$1.02$}} & \\
{\small{$\delta$}} & {\small{$1.03$}} &  &  &  & \\
\hline
\hline
\end{tabular}
\caption{Parameters of the TB-HIV/AIDS model \eqref{model:TB:HIV}.}
\label{table:parameters:TB-HIV}
\end{table}


\subsection{Equilibrium points and stability analysis}

In Table~\ref{table:effect:beta1} we show the effect
of the transmission coefficient $\beta_1$ on the state
$I_T^\diamond$ of the HIV-free equilibrium $\Sigma_T$
and on the basic reproduction number $R_1$.
Table~\ref{table:effect:beta2} shows the effect
of the transmission coefficient $\beta_2$ on the states
$I_H^\star$ and $A^\star$ of the TB-free equilibrium
$\Sigma_H$ and on the basic reproduction number $R_2$.
We conclude that the equilibrium states $I_T^\diamond$
and $(I_H^\star, A^\star)$ increase with the transmission
coefficients $\beta_1$ and $\beta_2$, respectively.
\begin{table}[!htb]
\centering
\begin{tabular}{ l   l   l   l   l  l}
\hline
{\small{$\beta_1$}}  &  {\small{4.3}} & {\small{6}}
&  {\small{10}}  &  {\small{15}} & {\small{50}}\\
\hline
\hline
{\small{$R_1$}}   &  {\small{0.99788}} & {\small{1.39239}}
&  {\small{2.32065}}  &  {\small{3.48097}} & {\small{11.60326}}\\
{\small{$I_T^\diamond$}}  &  {\small{0.00397}} & {\small{903.93492}}
&  {\small{2206.57268}}  &  {\small{2870.72755}} & {\small{3804.50589}}\\ \hline
\end{tabular}
\caption{Effect of $\beta_1$ on $I_T^\diamond$ and $R_1$.}
\label{table:effect:beta1}
\end{table}
\begin{table}[!htb]
\centering
\begin{tabular}{ l   l   l   l   l  l}
\hline
{\small{$\beta_2$}}  &  {\small{0.051}} & {\small{0.055}}
&  {\small{0.07}}  &  {\small{0.09}} & {\small{0.99}}\\
\hline
\hline
{\small{$R_2$}}   &  {\small{0.93669}} & {\small{1.01016}}
&  {\small{1.28566}}  &  {\small{1.65299}} & {\small{1.81829}}\\
{\small{$I_H^\star$}}  &  {\small{0.01708}} & {\small{135.73817}}
&  {\small{2516.54721}}  &  {\small{4472.84980}} & {\small{4930.48696}}\\
{\small{$A^\star$}}  &  {\small{0.00266}} & {\small{21.07182}}
&  {\small{390.59491}}  &  {\small{694.23361}} & {\small{765.26396}}\\
\hline
\end{tabular}
\caption{Effect of $\beta_2$ on $I_H^\star$, $A^\star$ and $R_2$.}
\label{table:effect:beta2}
\end{table}


\bigskip

\begin{figure}[!htb]
$\begin{array}{cc}
\hspace*{-0.4in}
\includegraphics[scale=0.6]{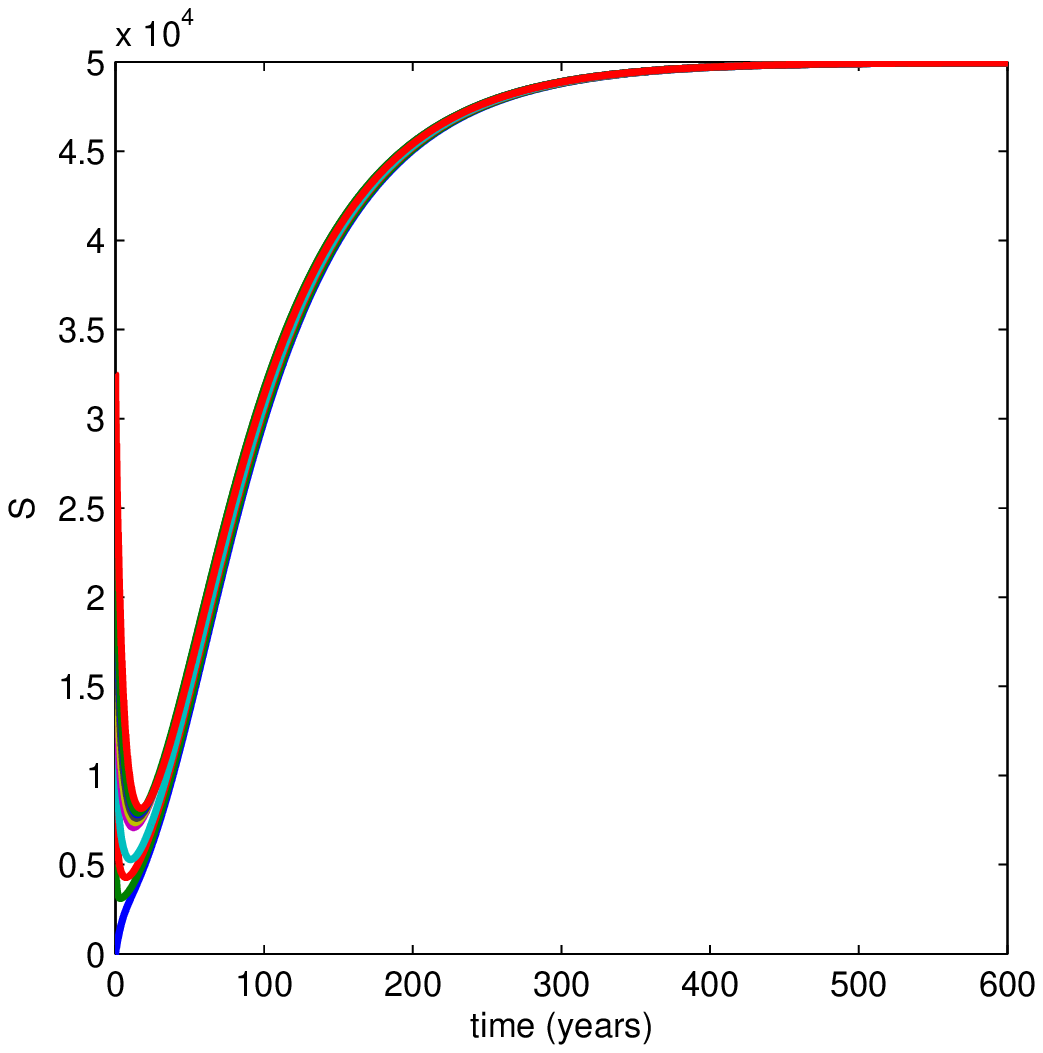} &
\hspace*{-0.6in}\includegraphics[scale=0.6]{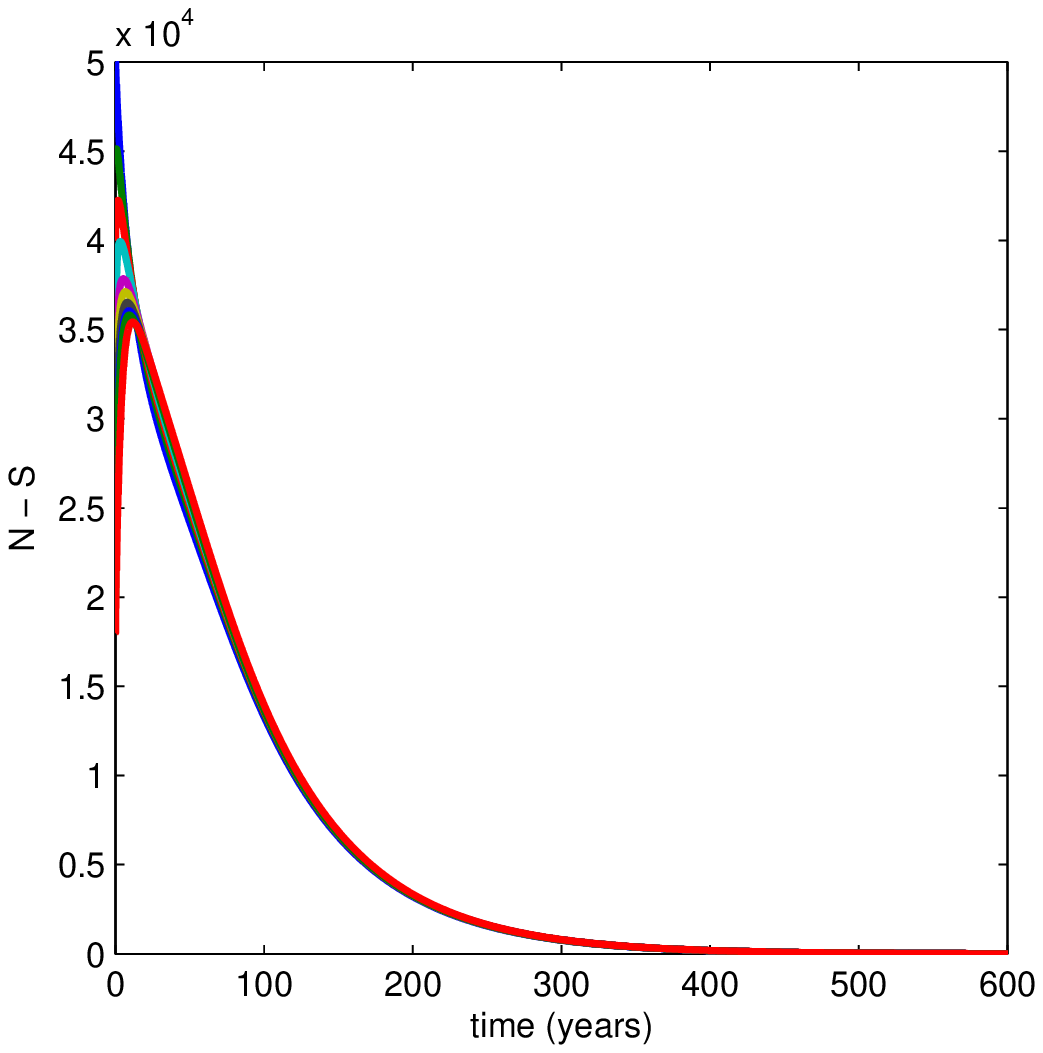}\\
\end{array}$
\caption{Stability of the disease-free equilibrium \eqref{eq:DFE:model:TBHIV}.}
\label{fig:DFE:stab}
\end{figure}

In Figure~\ref{fig:DFE:stab} we considered different initial conditions
in a neighborhood of the initial conditions given by \eqref{init:cond}
and $R_0 < 1$ ($R_1 < 1$ and $R_2 < 1$) to illustrate the stability
of the disease-free equilibrium $\Sigma_0$ given by \eqref{eq:DFE:model:TBHIV}.
In these numerical simulations we considered $\beta_1 = 2.7$ and $\beta_2 = 0.03$,
corresponding to $R_1 = 0.62632$ and $R_2 = 0.55077$, while the rest of
the parameters take the values in Table~\ref{table:parameters:TB-HIV}.

\begin{figure}[!htb]
\hspace*{-0.6in}
  \includegraphics[scale=0.39]{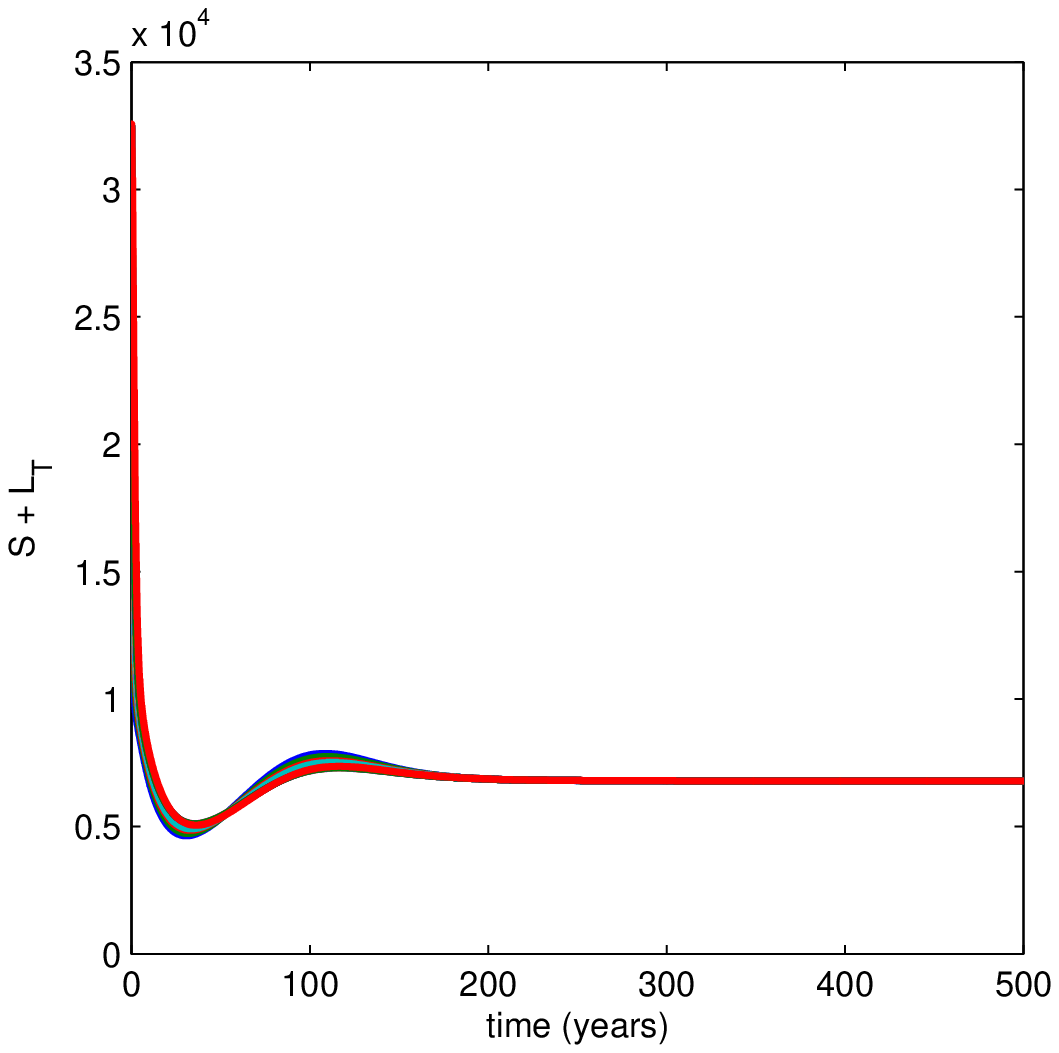}
  \hspace*{-1.1in}
  \includegraphics[scale=0.39]{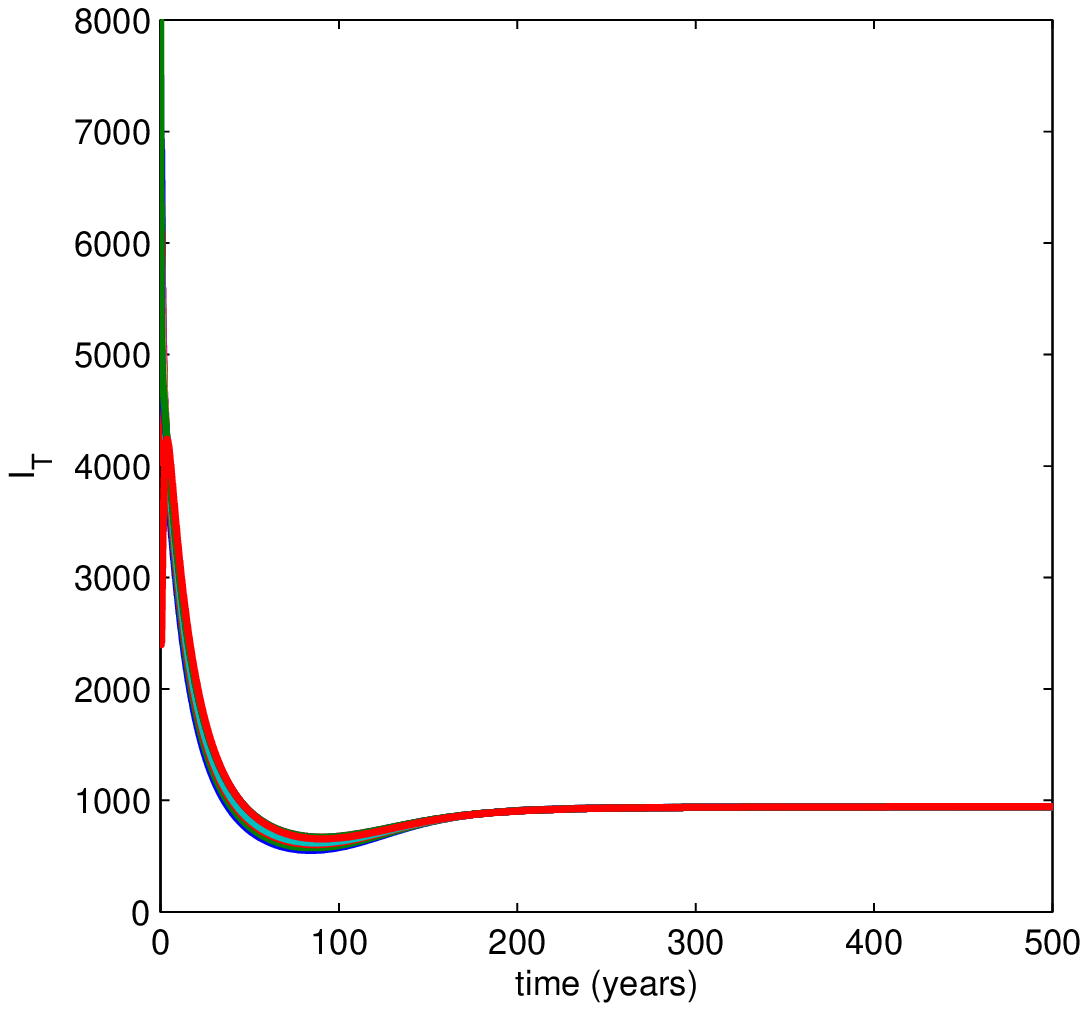}
  \hspace*{-1.1in}
  \includegraphics[scale=0.39]{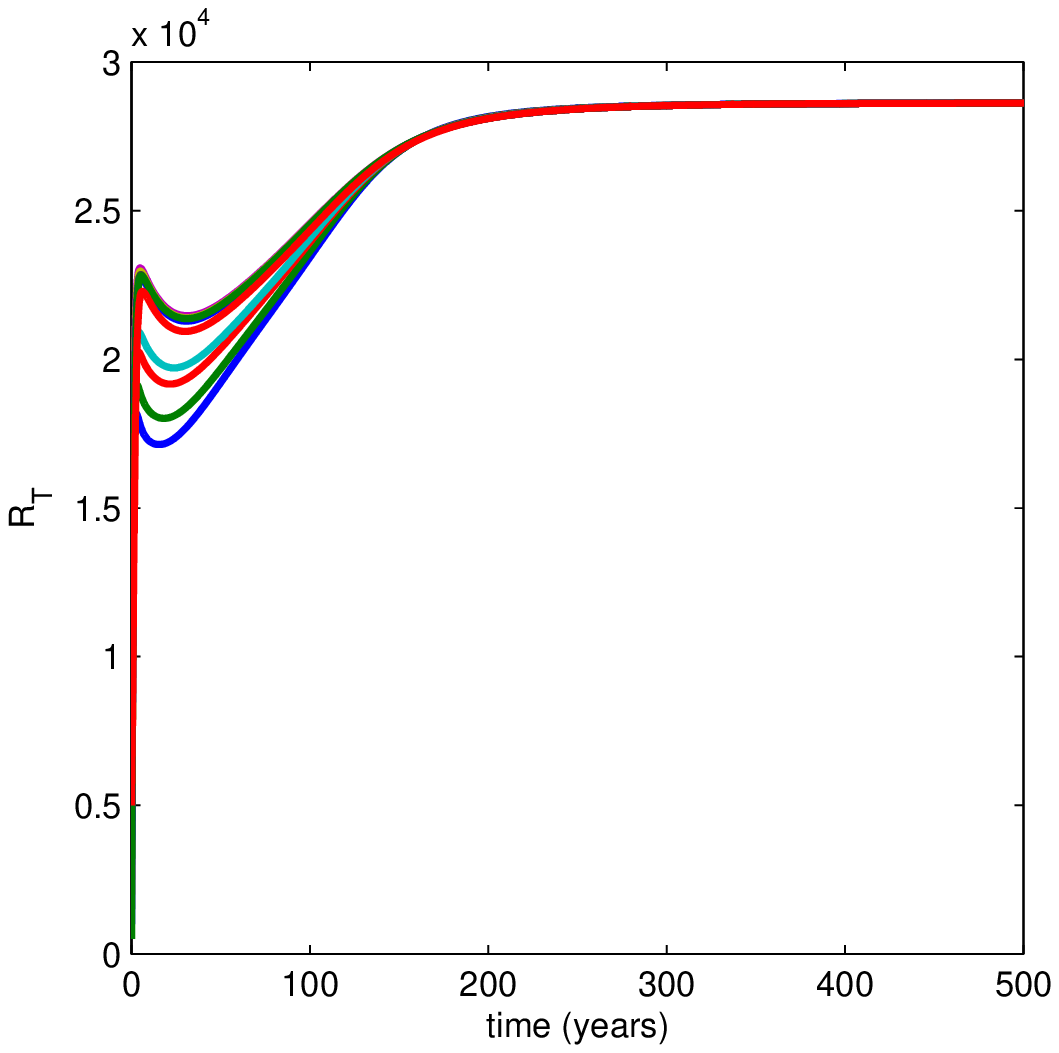}
  \hspace*{-0.28in}
  \includegraphics[scale=0.39]{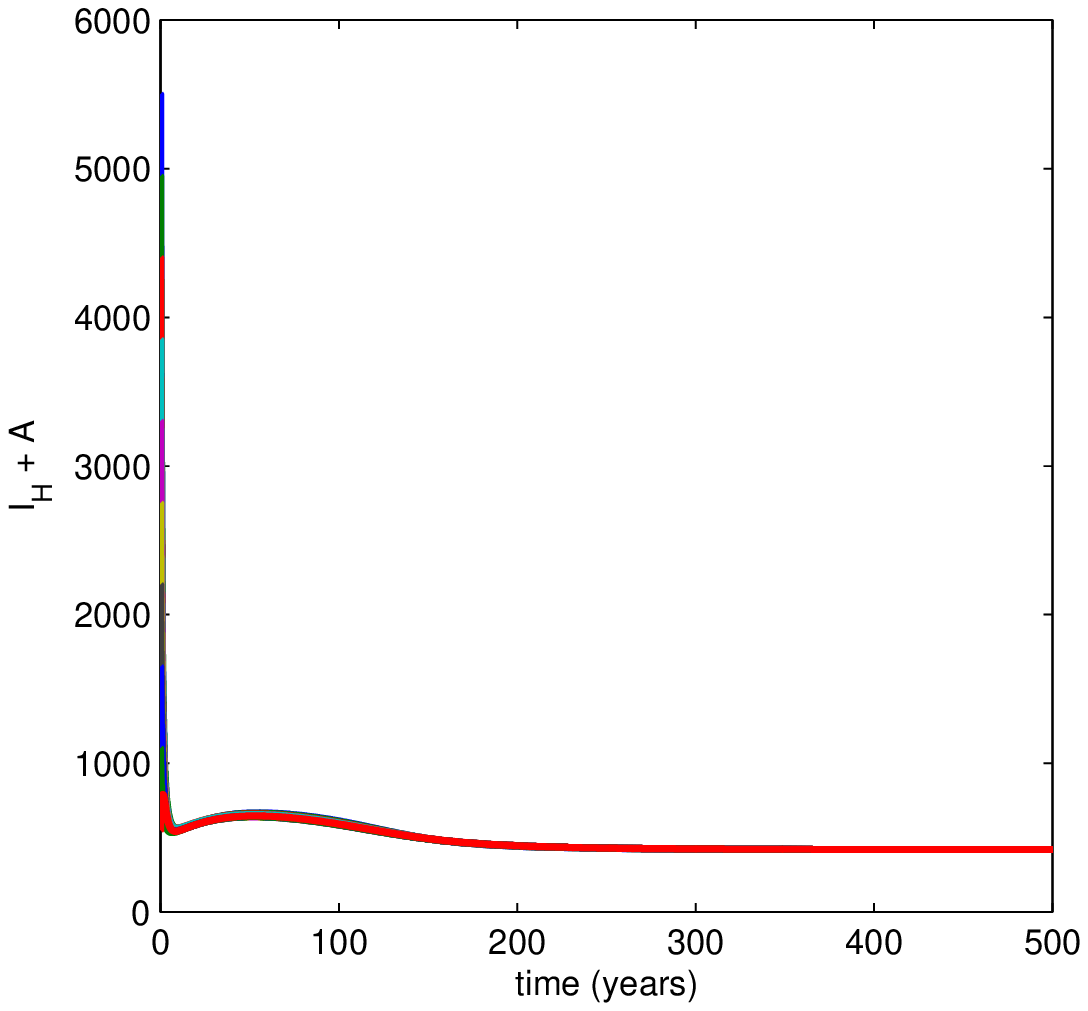}
  \hspace*{-0.45in}
  \includegraphics[scale=0.39]{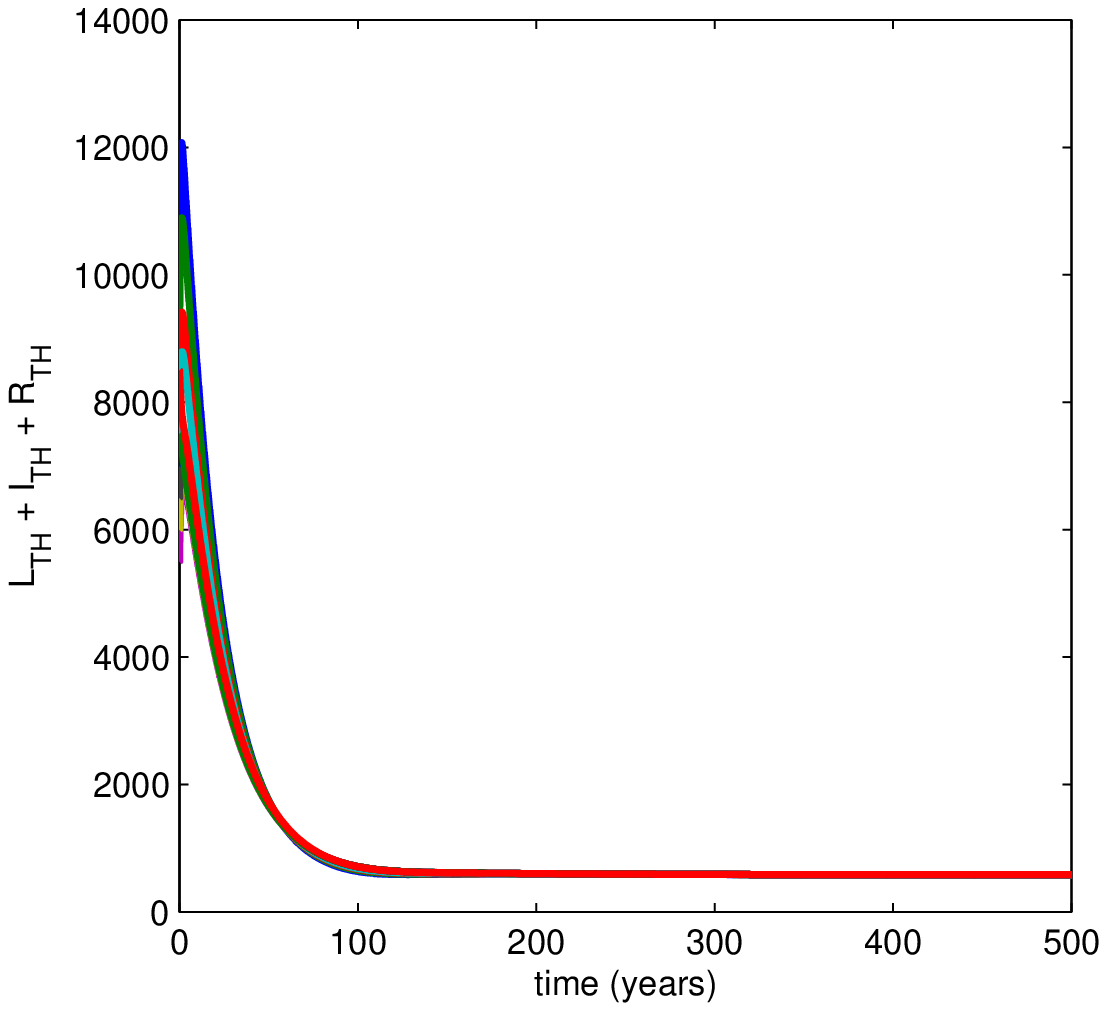}
  \hspace*{-0.45in}
  \includegraphics[scale=0.39]{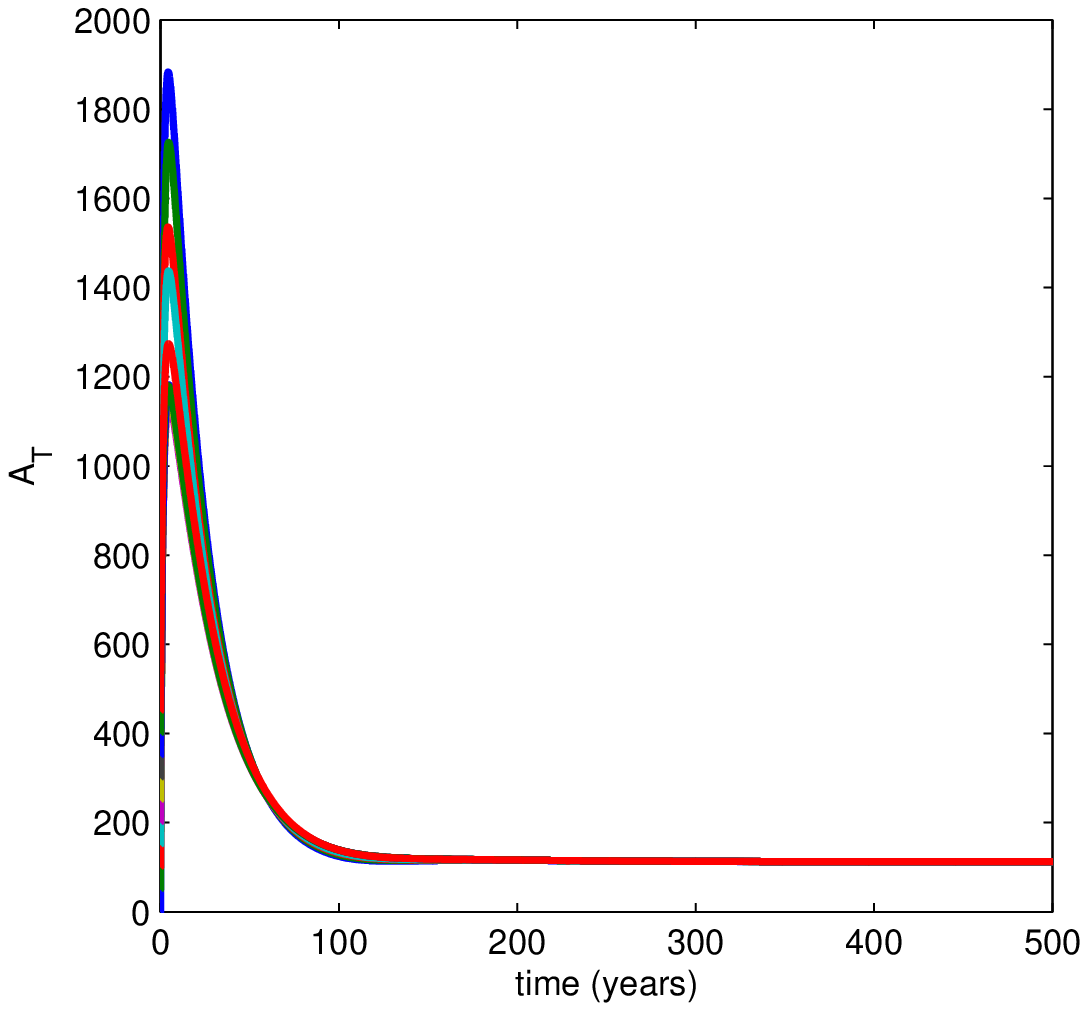}
\caption{Stability of the syndemic equilibrium $\Sigma^*$.}
\label{fig:EE:stab}
\end{figure}

Figure~\ref{fig:EE:stab} shows that, for $R_0 > 1$, the syndemic equilibrium
$\Sigma^*$ exists. We considered different initial conditions for the state
variables of system \eqref{model:TB:HIV} in a neighborhood of \eqref{init:cond},
$\beta_1 = 6$ and $\beta_2 = 0.1$, corresponding to $R_1 = 1.39239$
and $R_2 = 1.83593$, and the rest of the parameters take the values
in Table~\ref{table:parameters:TB-HIV}. We observe that the state variables
converge to $\Sigma^*$ when $t \rightarrow \infty$. In this case, $\Sigma^*$ is given by
\begin{equation*}
\begin{split}
\Sigma^* &= \left(S^*, L_T^*, I_T^*, R_T^*, I_H^*, A^*, L^*_{TH}, I_{TH}^*, R_{TH}^*, A_T^*\right)\\
&= (4766.84, 2019.66, 943.06, 28621.89, 362.66, 56.29, 31.39, 55.15, 495.68, 112.33) \, .
\end{split}
\end{equation*}

\subsection{Treatment impact on TB-HIV/AIDS coinfection}

\begin{figure}[!htb]
\hspace*{-0.22in}
  \includegraphics[scale=0.39]{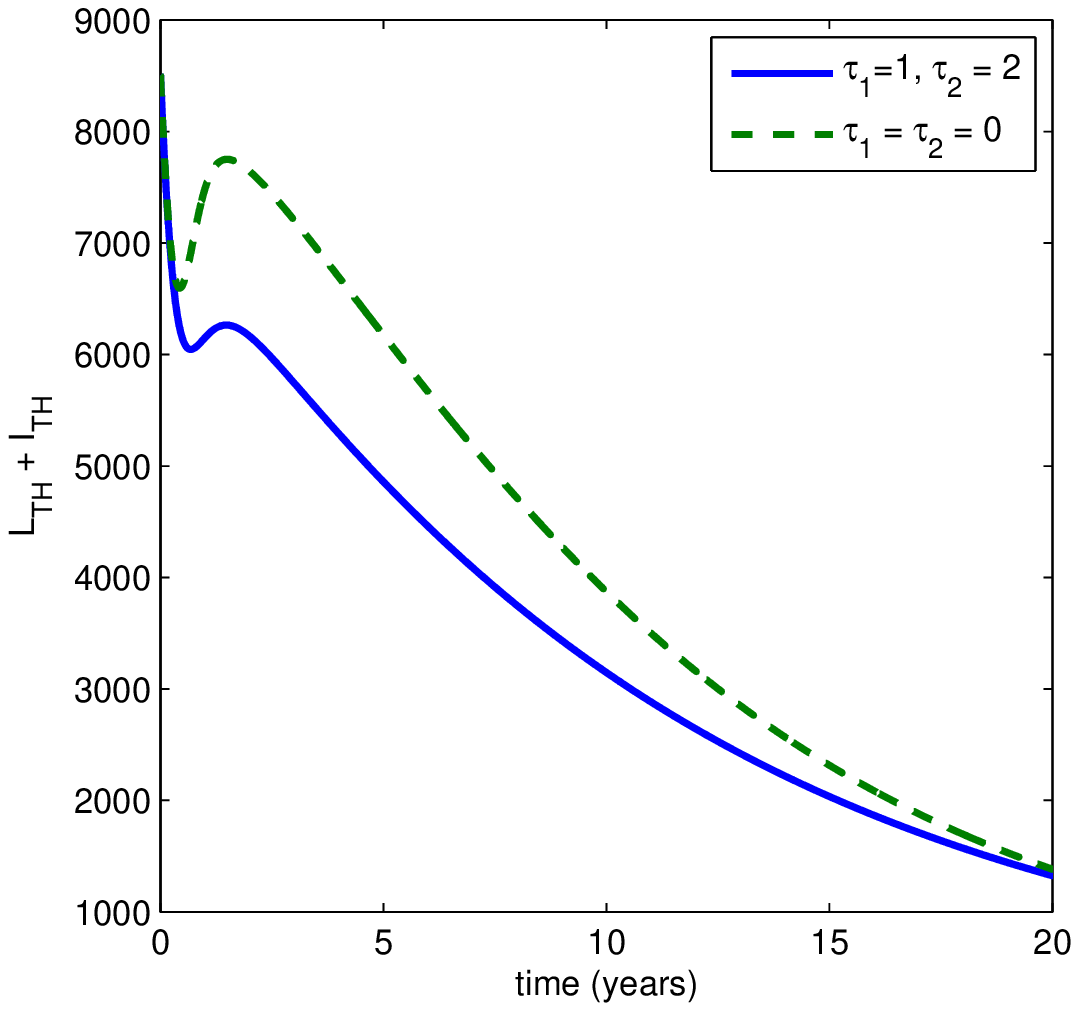}
  \hspace*{-0.45in}
  \includegraphics[scale=0.39]{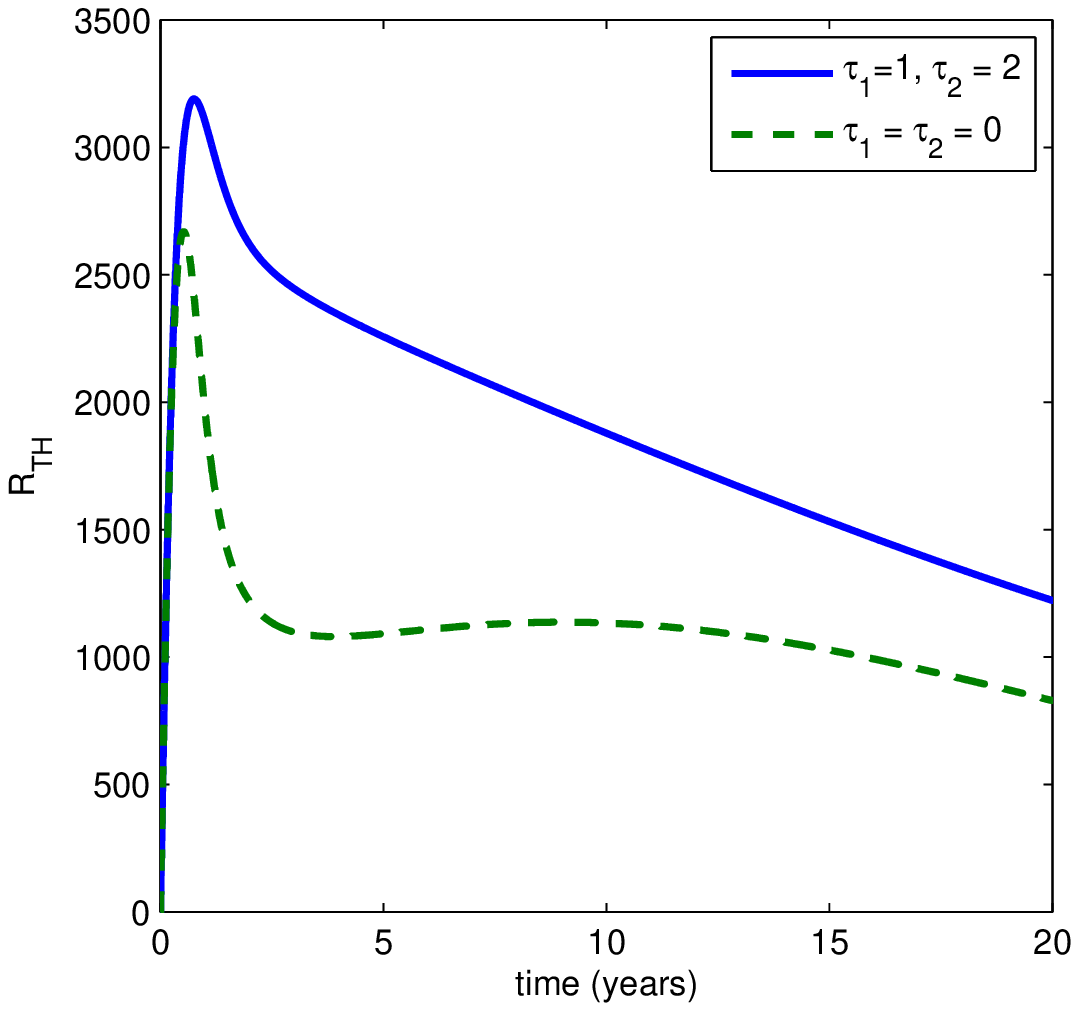}
  \hspace*{-0.45in}
  \includegraphics[scale=0.39]{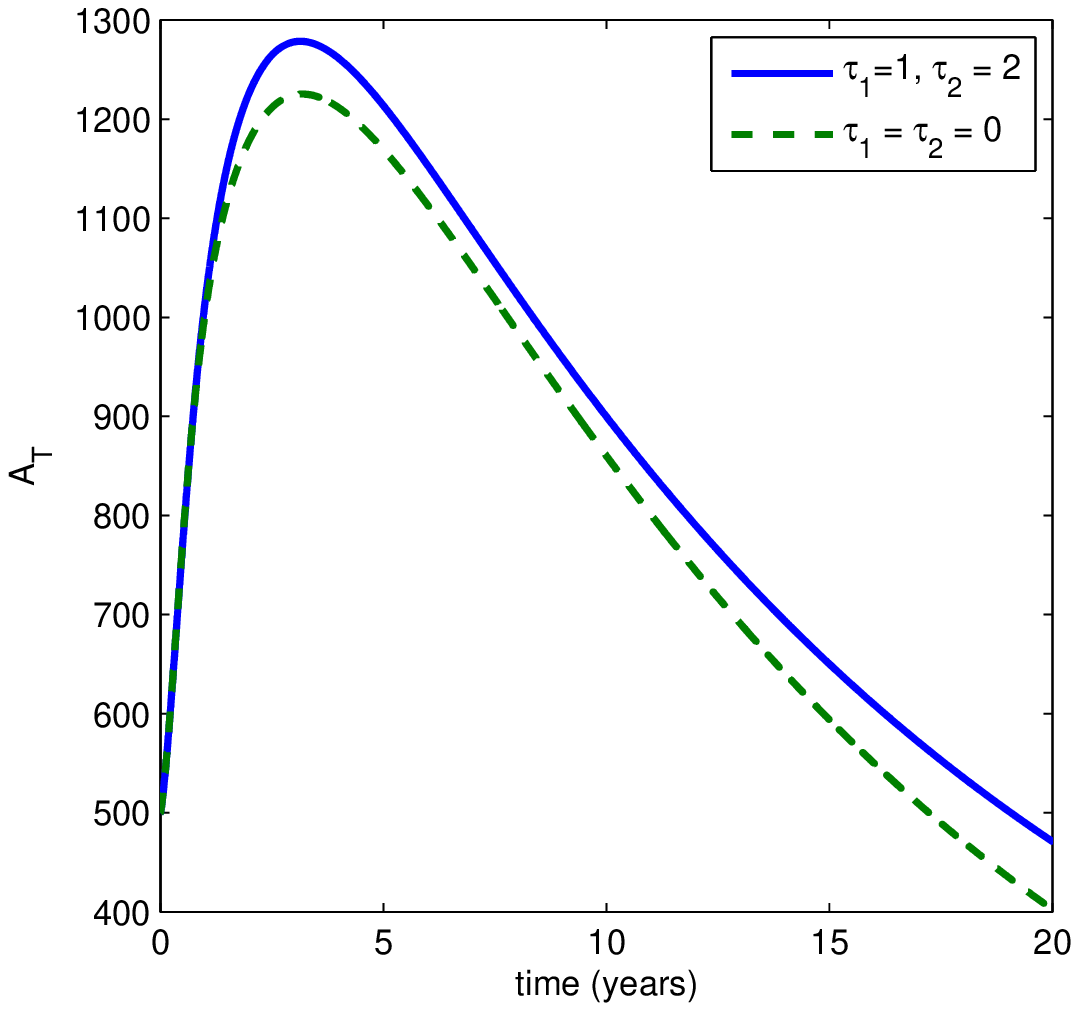}
\caption{Impact of TB treatment on single-infected individuals with disease induced death.}
\label{fig:treatTBonly:death}
\end{figure}

\begin{figure}[!htb]
\hspace*{-0.22in}
  \includegraphics[scale=0.39]{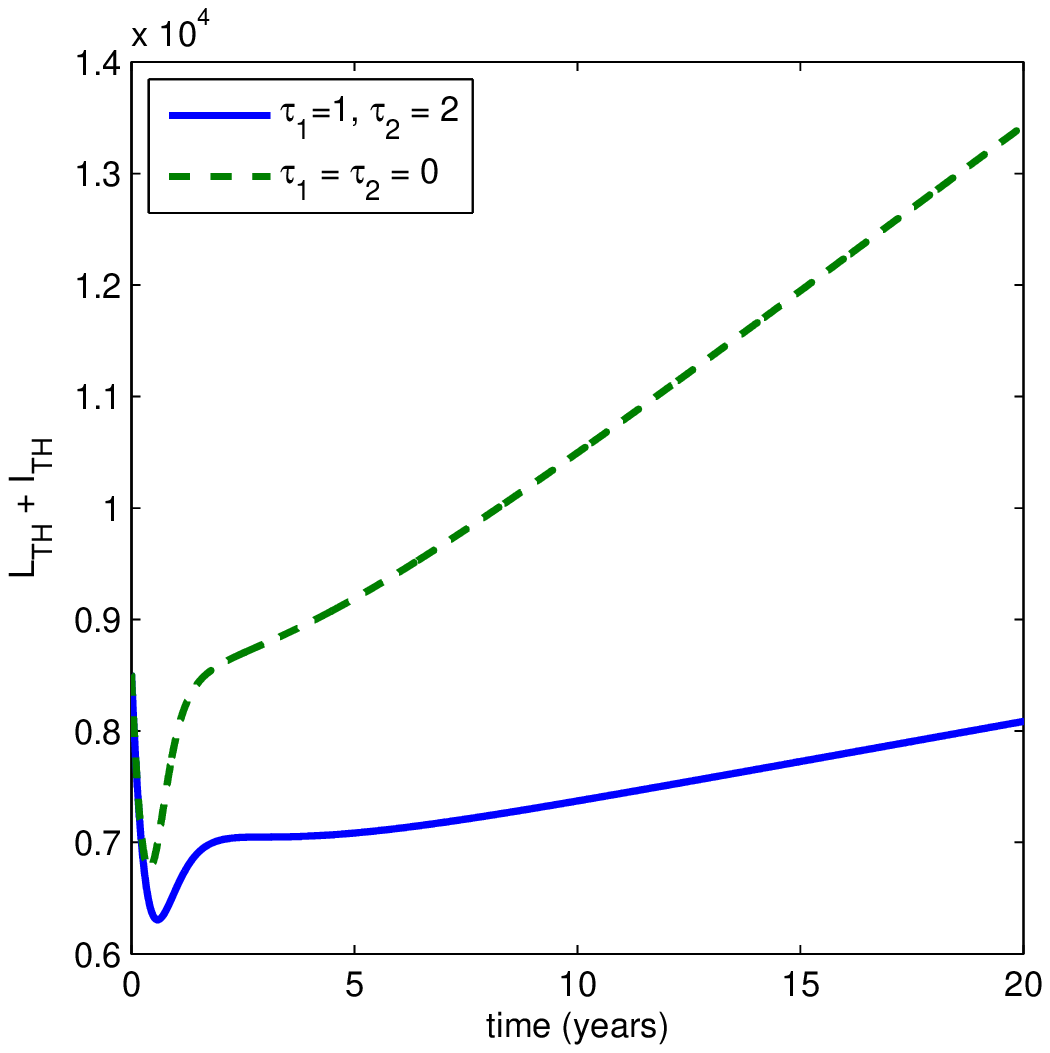}
  \hspace*{-0.45in}
  \includegraphics[scale=0.39]{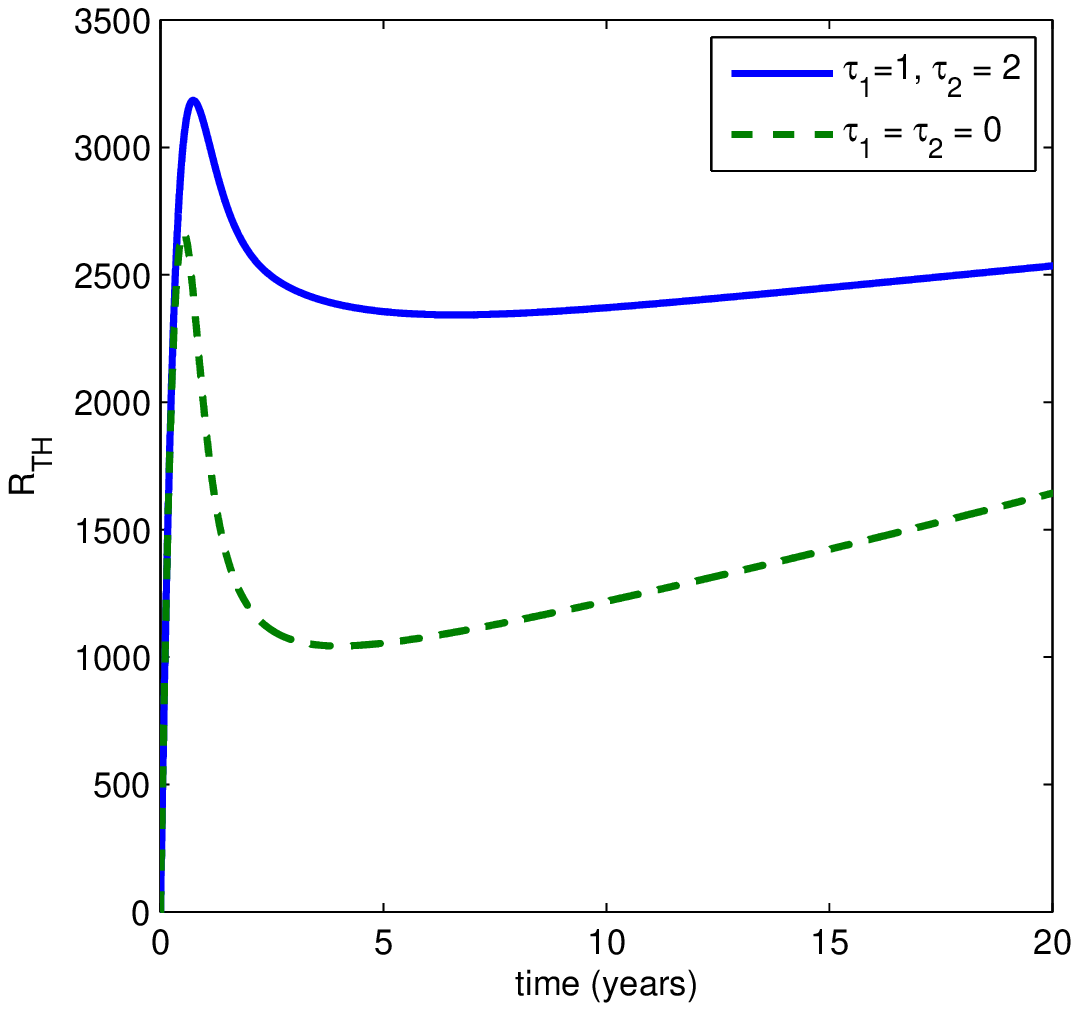}
  \hspace*{-0.45in}
  \includegraphics[scale=0.39]{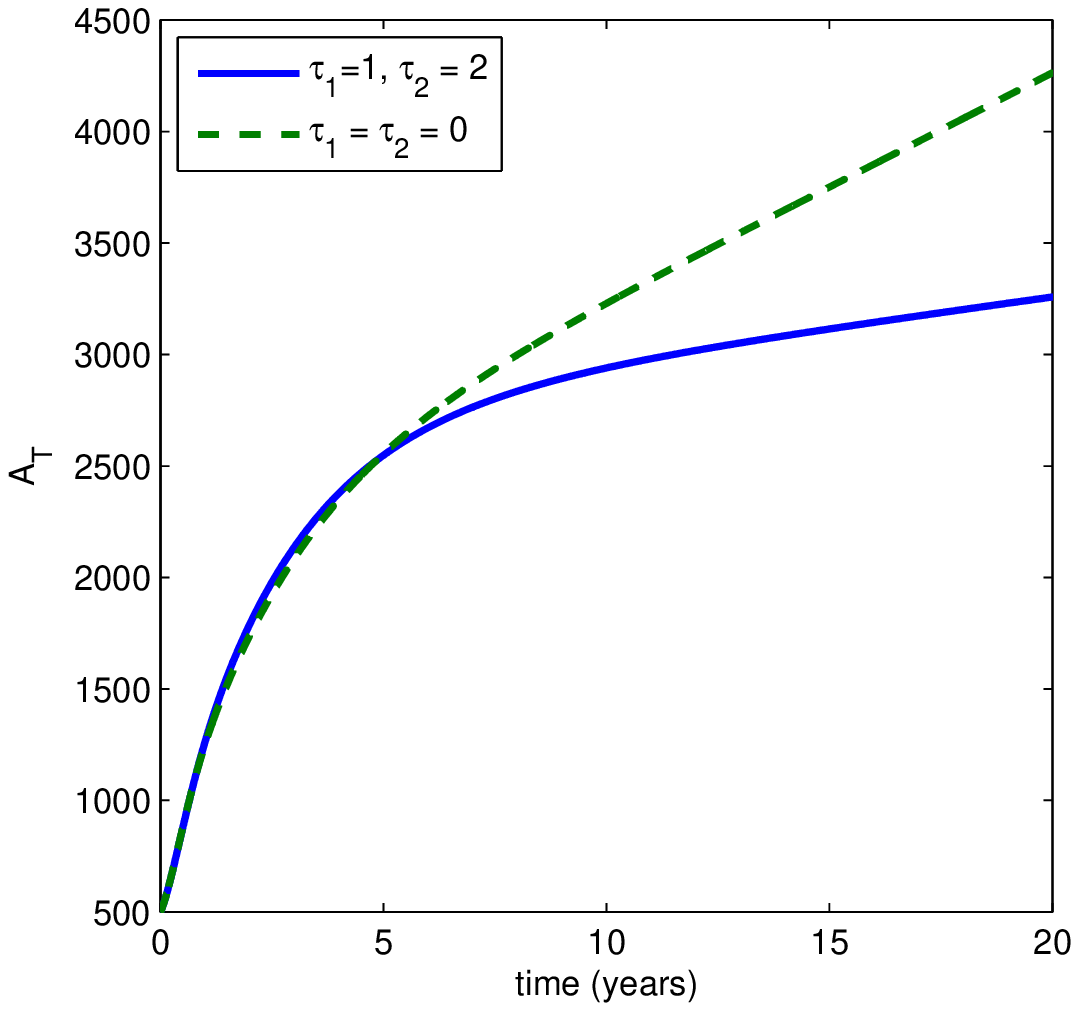}
\caption{Impact of TB treatment on single-infected individuals with no disease induced death.}
\label{fig:treatTBonly:no:death}
\end{figure}

\begin{figure}[!htb]
\hspace*{-0.22in}
  \includegraphics[scale=0.39]{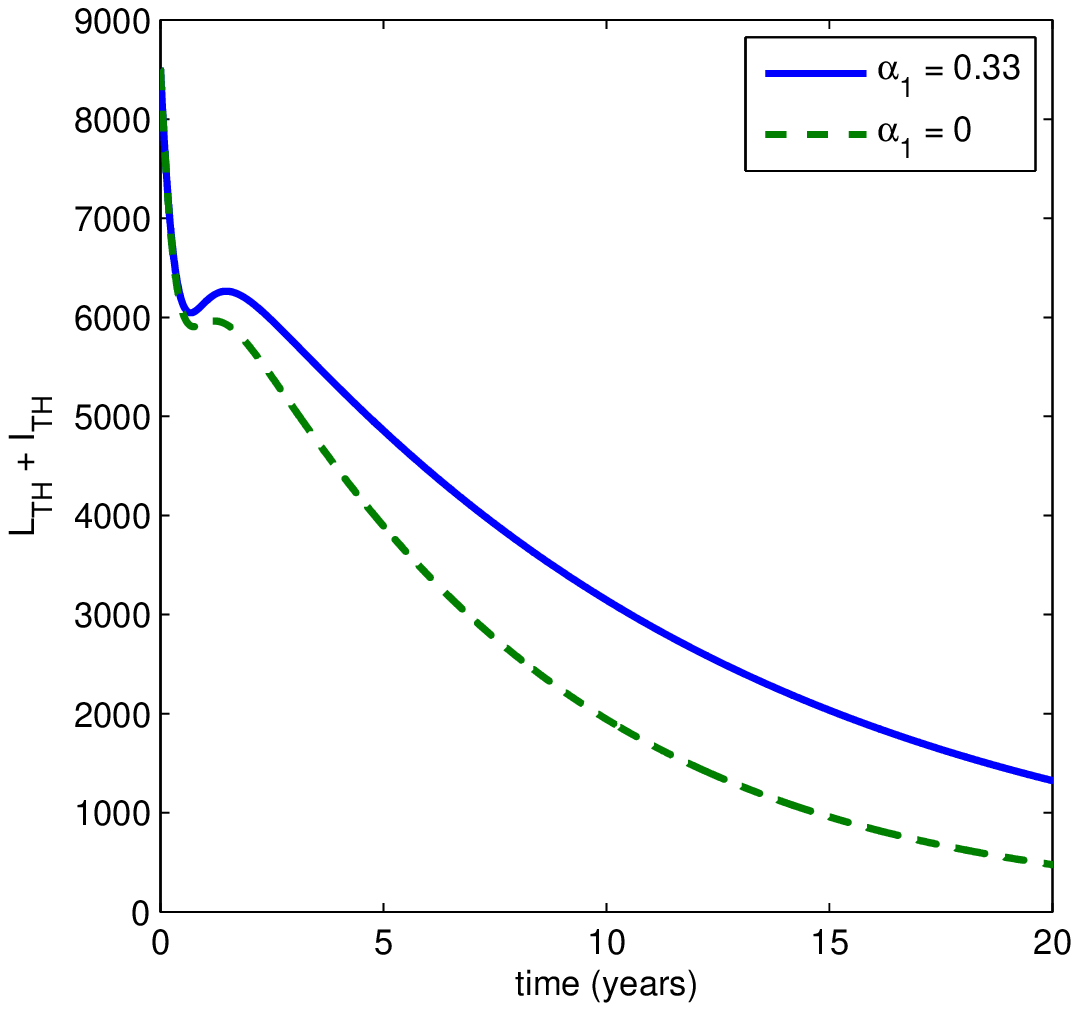}
  \hspace*{-0.45in}
  \includegraphics[scale=0.39]{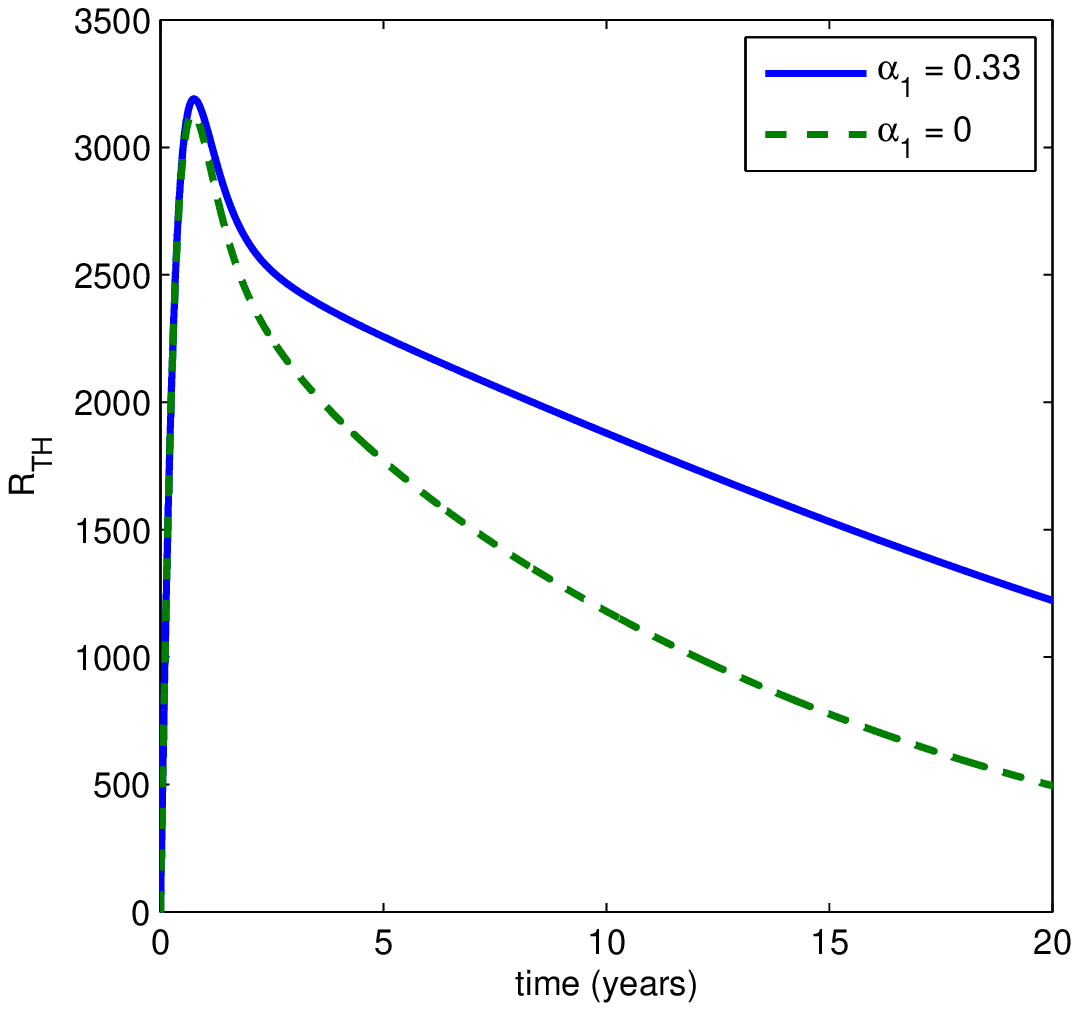}
  \hspace*{-0.45in}
  \includegraphics[scale=0.39]{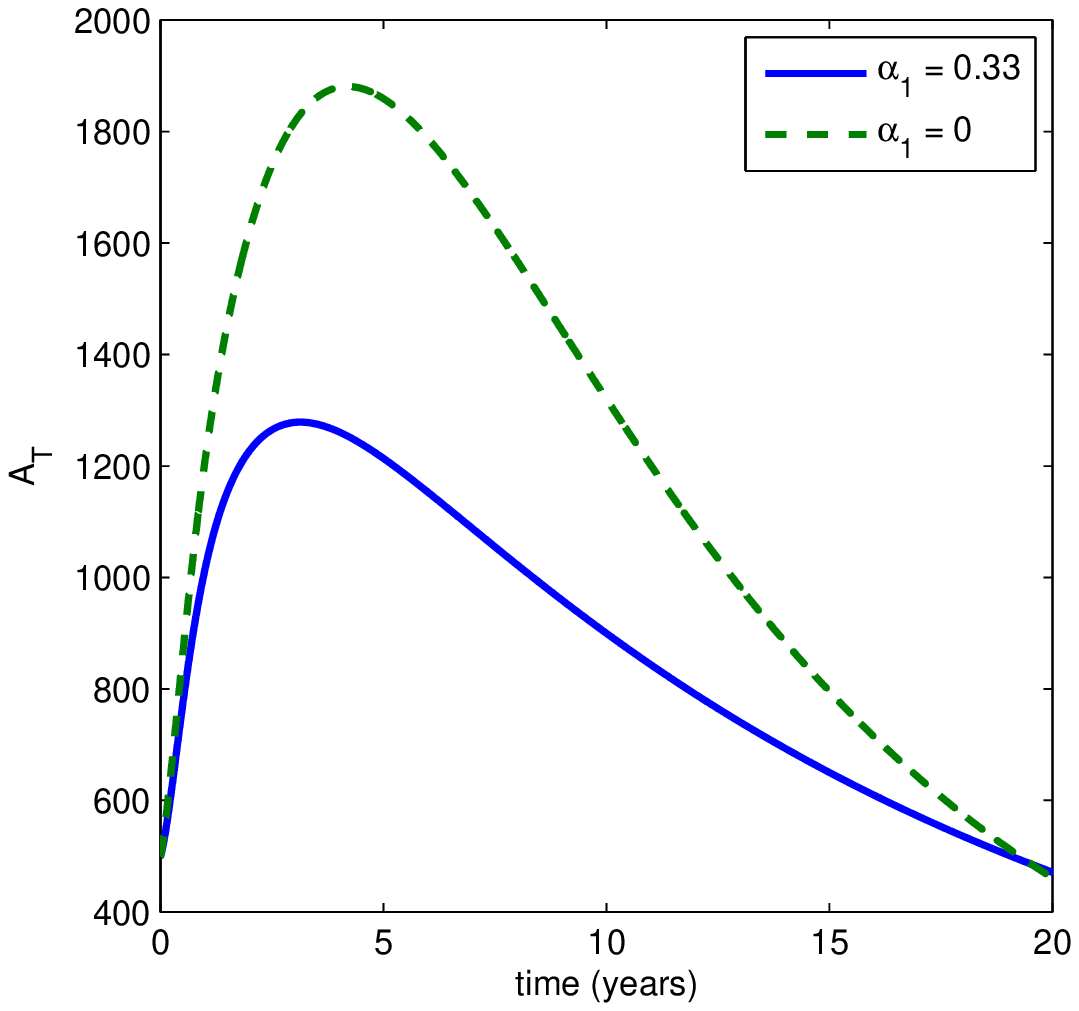}
\caption{Impact of AIDS treatment on single-infected individuals with disease induced death.}
\label{fig:treatHIVonly:death}
\end{figure}

\begin{figure}[!htb]
\hspace*{-0.22in}
  \includegraphics[scale=0.39]{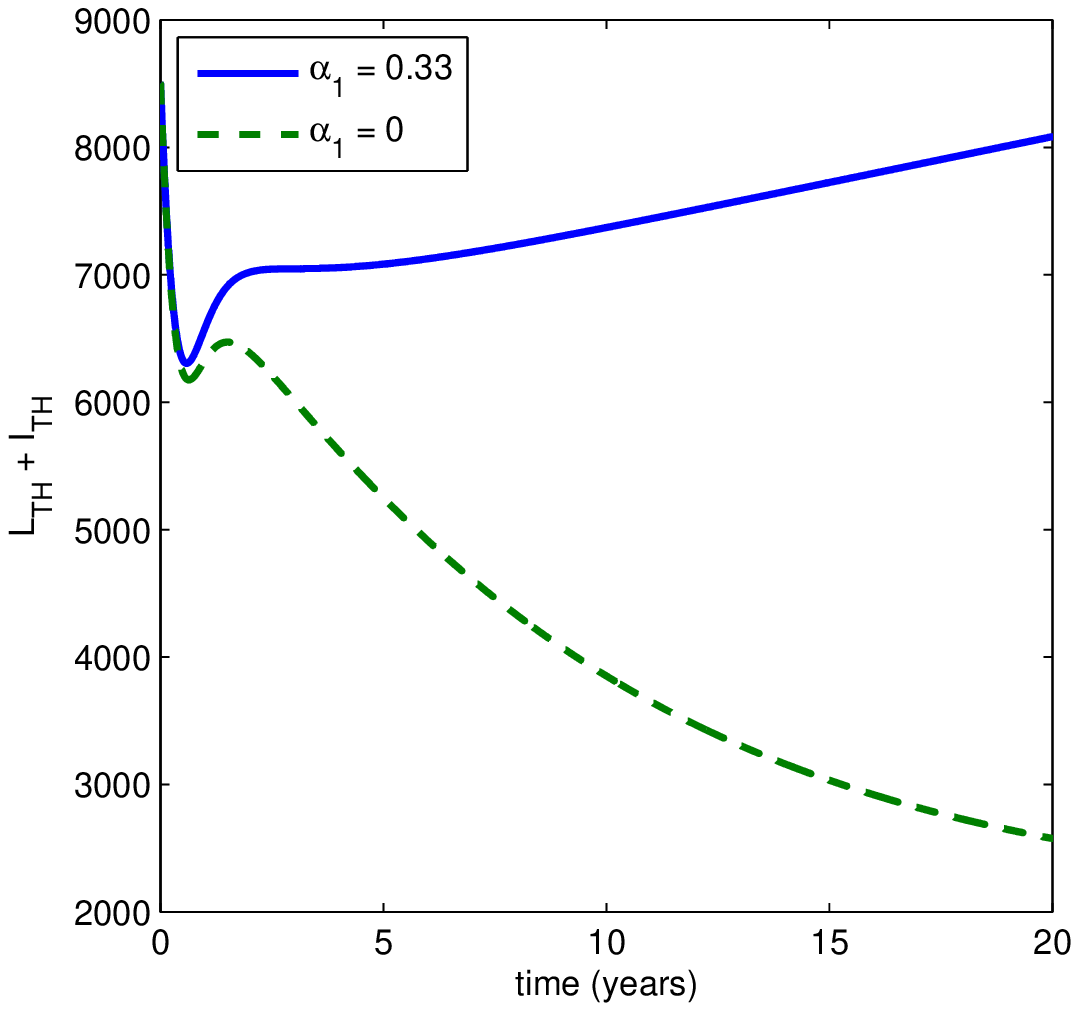}
  \hspace*{-0.45in}
  \includegraphics[scale=0.39]{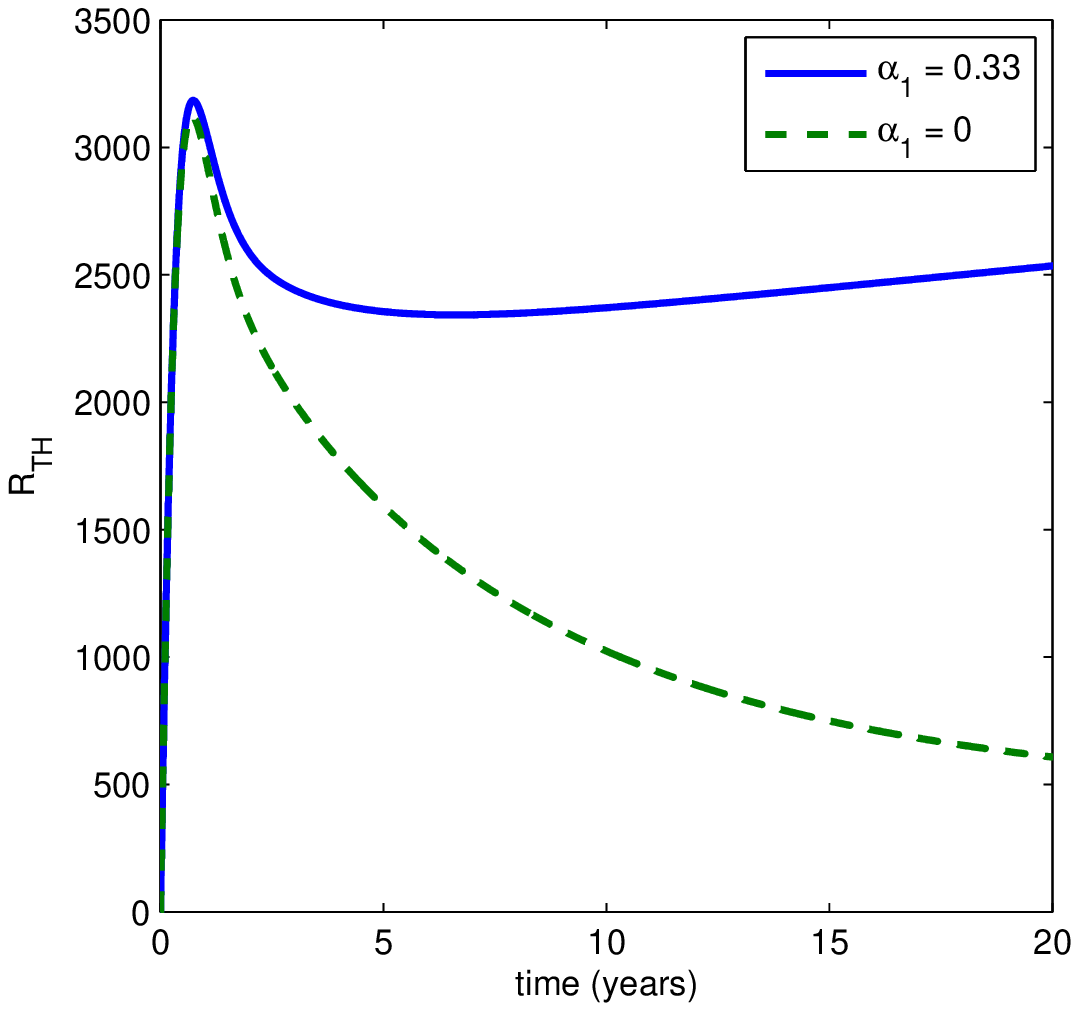}
  \hspace*{-0.45in}
  \includegraphics[scale=0.39]{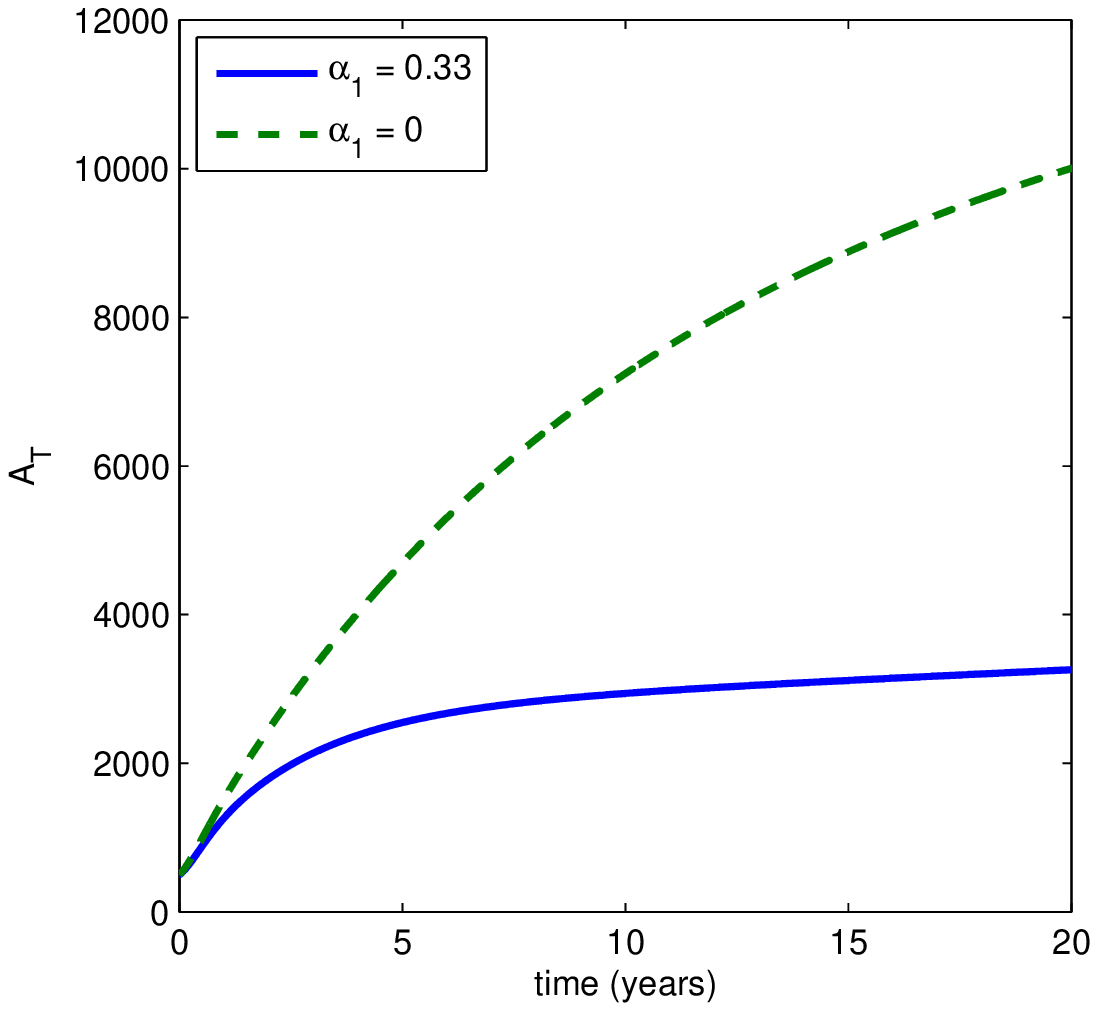}
\caption{Impact of AIDS treatment on single-infected individuals with no disease induced death.}
\label{fig:treatHIVonly:no:death}
\end{figure}

\begin{figure}[!htb]
\hspace*{-0.22in}
  \includegraphics[scale=0.39]{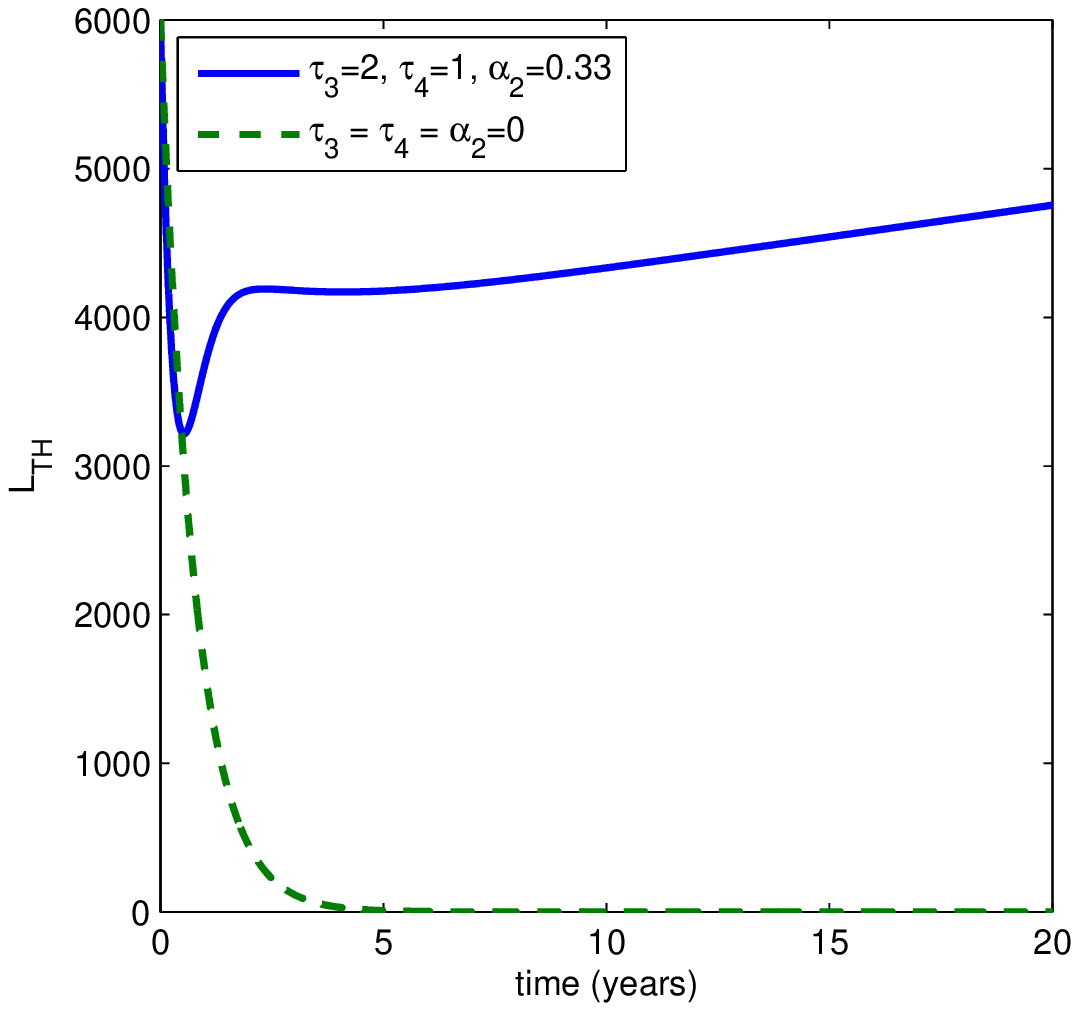}
  \hspace*{-0.45in}
  \includegraphics[scale=0.39]{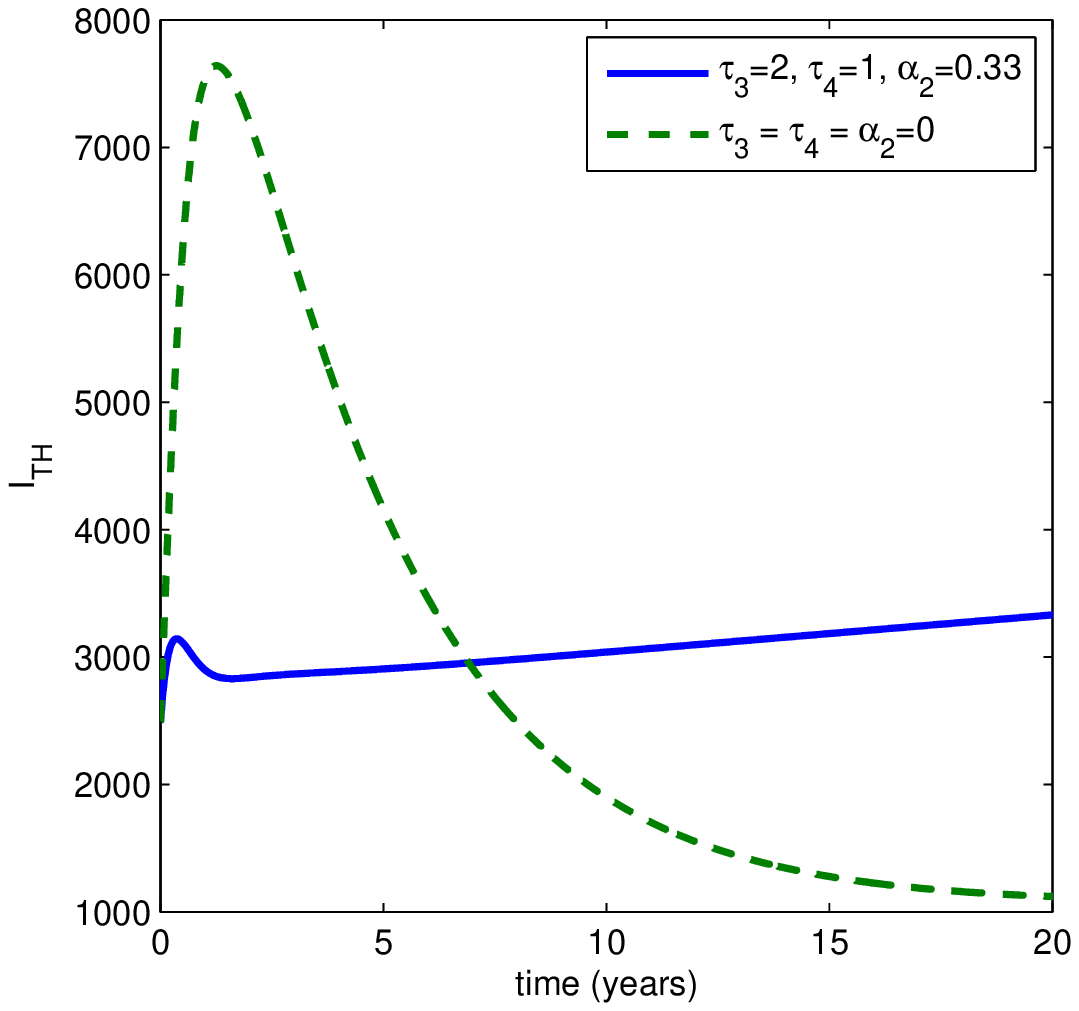}
  \hspace*{-0.45in}
  \includegraphics[scale=0.39]{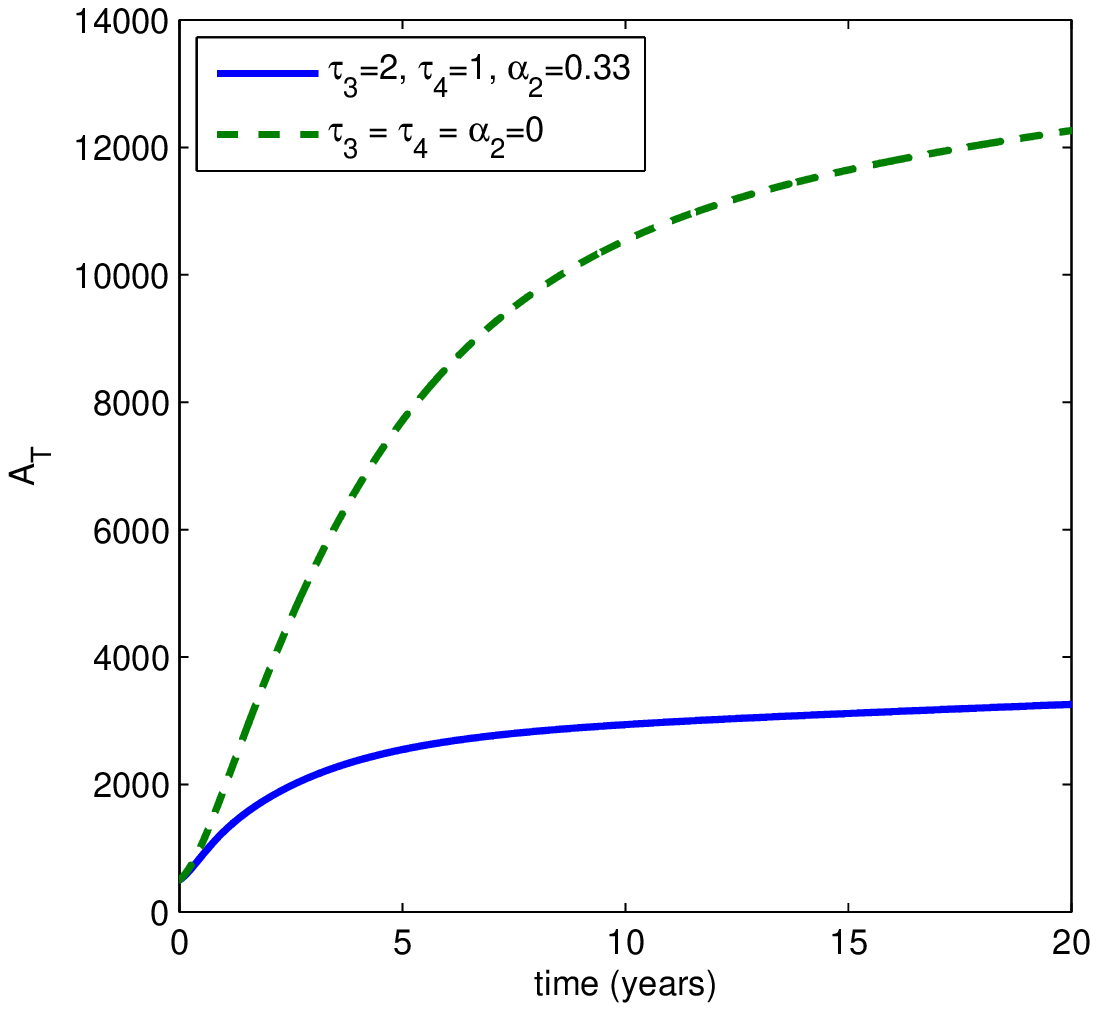}
\caption{Impact of TB and AIDS treatment on co-infected individuals with no disease induced death.}
\label{fig:treatTBHIV:no:death}
\end{figure}

Consider $\beta_1=13$ and $\beta_2 = 0.06$, while the rest of the parameters
take the values of Table~\ref{table:parameters:TB-HIV}.
Figure~\ref{fig:treatTBonly:death} shows the impact of treating
the individuals with active and latent TB on the number of individuals
co-infected with TB-HIV/AIDS. The treatment of individuals with only-TB,
$I_T$ and $L_T$, has a positive impact on the reduction of the number
of individuals co-infected with TB-HIV/AIDS. Moreover, the number
of individuals that suffered from disease (TB and AIDS) induced death
is higher when individuals with TB-single infection are not treated.
In this case the total population at the end of 20 years is around
10509 and, in the case where individuals with only TB are treated,
the total population at the end of 20 years is around 29758 individuals.
In Figure~\ref{fig:treatTBonly:no:death}, we assume that there are no
disease induced deaths, that is, $d_T = d_A = d_{TA} = 0$. The impact
of treating individuals with only TB on the reduction of the number
of individuals co-infected is more evident.
Figure~\ref{fig:treatHIVonly:death} illustrates the case where we compare the number
of individuals co-infected with TB-HIV/AIDS when individuals with only AIDS symptoms
$A_T$ are or not treated. We observe that treating this class of individuals is important
for the reduction of the number of individuals that become co-infected, with special
attention to the individuals that have AIDS symptoms and TB infection.
In figure~\ref{fig:treatHIVonly:no:death}, we considered that there
is no disease induced deaths ($d_T = d_A = d_{TA} = 0$).
It is crucial that TB-infected individuals (in the latent and active stage), which
are also HIV-positive, take anti-TB drugs, since they can recover from TB.
We analyze the impact of treating TB-HIV/AIDS co-infected individuals
$L_{TH}$, $I_{TH}$ and $A_T$ on the reduction of the number of individuals
coinfection. If anti-TB drugs are supplied, then latent and active-TB individuals
with HIV can recover and pass to the class $R_{TH}$ (the number of individuals
in the class $R_{TH}$ tends to zero when TB is not treated).
In Figure~\ref{fig:treatTBHIV:no:death}, we observe that after 7 years
the number of individuals infected with active-TB and HIV, in the case
without treatment, becomes lower than in the case with treatment.
This is due to the fact that coinfection precipitates AIDS symptoms.


\appendix


\section{Appendix: Computation of $R_0$}
\label{A.1}

The basic reproduction number represents the expected average number of new infections
produced by a single infectious individual when in contact with a completely susceptible
population \cite{van:den:Driessche:2002}.
Following \cite{van:den:Driessche:2002}, the basic reproduction number $R_{0}$ is obtained
as the spectral radius of the matrix $F \cdot V^{-1}$ at the disease-free equilibrium
$\Sigma_0$, given by \eqref{eq:DFE:model:TBHIV}, with
$F = \left[ \begin{array}{cc} F_1 \quad F_2 \end{array} \right]$ and
{\footnotesize{
\begin{equation*}
F_1 =
 \left[ \begin {array}{ccccc}
0&0&0&0&0\\ \noalign{\medskip}
\lambda_T&0&{\frac {\beta_1 S}{N}}+{\frac{\beta^{'}_1 \beta_1 R_T}{N}}&
\beta^{'}_1 \lambda_T&0\\ \noalign{\medskip}
0&0&0&0&0\\ \noalign{\medskip}
0&0&0&0&0\\ \noalign{\medskip}
\lambda_H&0&0&\lambda_H&{\frac{\beta_2 S}{N}}+{\frac{\beta_2 R_T}{N}}\\ \noalign{\medskip}
0&0&0&0&0\\ \noalign{\medskip}
0&0&{\frac{\beta^{'}_2 \beta_1 R_{TH}}{N}}&0&0\\ \noalign{\medskip}
0&0&\delta \lambda_H+{\frac{\psi \beta_1 I_H}{N}}&0&
{\frac{\delta \beta_2 I_T}{N}}+ \psi \lambda_T\\ \noalign{\medskip}
0&0&0&0&0\\ \noalign{\medskip}
0&0&0&0&0\end {array} \right]
\end{equation*}
}}
{\footnotesize{
\begin{equation*}
F_2 =
 \left[ \begin {array}{cccccccccc}
0&0&0&0&0\\ \noalign{\medskip}
0&0&{\frac{\beta_1 S}{N}}+{\frac{\beta^{'}_1 \beta_1 R_T}{N}}
&0&\frac{\beta_1 S}{N} + \frac{\beta_1^{'} \beta_1 R_T}{N}\\
\noalign{\medskip}
0&0&0&0&0\\ \noalign{\medskip}
0&0&0&0&0\\ \noalign{\medskip}
{\frac{\beta_2 \eta S}{N}}+{\frac{\beta_2 \eta R_T}{N}}&
{\frac{\beta_2 S}{N}}+{\frac{\beta_2 R_T}{N}}&
{\frac{\beta_2 S}{N}}+{\frac{\beta_2 R_T}{N}}&
{\frac{\beta_2 S}{N}}+{\frac{\beta_2 R_T}{N}}&
{\frac{\beta_2 S \eta }{N}}+{\frac{R_T \beta_2 \eta }{N}}\\
\noalign{\medskip}
0&0&0&0&0\\ \noalign{\medskip}
0&0&{\frac{\beta^{'}_2 \beta1 R_{TH}}{N}}& \beta^{'}_2 \lambda_T
&\frac{\beta_2^{'}\beta_1 R_{TH}}{N}\\
\noalign{\medskip}
{\frac {\delta \beta_2 \eta I_T}{N}}
&{\frac{\delta \beta_2 I_T}{N}}
&{\frac{\delta \beta_2 I_T}{N}}+{\frac {\psi \beta_1 I_H}{N}}
&{\frac{\delta \beta_2 I_T}{N}}
&{\frac {\delta \beta_2 \eta  I_T}{N}}+ \frac{\psi \beta_1 I_H}{N}\\
\noalign{\medskip}
0&0&0&0&0\\ \noalign{\medskip}
0&0&0&0&0\end {array} \right]
\end{equation*}
}}
and
$V = \left[ \begin{array}{cc} V_1 \quad V_2 \end{array} \right]$
with
{\footnotesize{
\begin{equation*}
V_1 =
\left[
\begin {array}{cccccccccc}
\lambda_T +\lambda_H+\mu&0&{\frac{\beta_1 S}{N}}&0&{
\frac{\beta_2 S}{N}}\\ \noalign{\medskip}
0&k_1 \tau_1+\mu&0&0&0\\ \noalign{\medskip}
0&-k_1 &\tau_2+\delta \lambda_H+\mu+d_T&0&{
\frac{\delta \beta_2 I_T}{N}}\\ \noalign{\medskip}
0&-\tau_1&-\tau_2+{\frac{\beta^{'}_1 \beta_1 R_T}{N}}&
\beta^{'}_1 \lambda_T+\lambda_H+\mu
&{\frac {\beta_2 R_T}{N}}\\ \noalign{\medskip}
0&0&{\frac {\psi \beta_1 I_H}{N}}&0&\rho_1+\psi \lambda_T +\mu\\
\noalign{\medskip}
0&0&0&0&-\rho_1\\ \noalign{\medskip}
0&0&0&0&0\\ \noalign{\medskip}
0&0&0&0&0\\ \noalign{\medskip}
0&0&{\frac{\beta^{'}_2 \beta_1 R_{TH}}{N}}&0&0\\
\noalign{\medskip}
0&0&0&0&0\end {array} \right]
\end{equation*}
}}
{\footnotesize{
\begin{equation*}
V_2 =
 \left[ \begin {array}{ccccc}
{\frac{\beta_2 \eta S}{N}}&
{\frac{\beta_2 S}{N}}&
{\frac{\beta_1 S}{N}}+{\frac{\beta_2 S}{N}}&
{\frac{\beta_2 S}{N}}
&\frac{\beta_1 S}{N} +{\frac{\beta_2 S \eta }{N}}\\
\noalign{\medskip}
0&0&0&0&0\\ \noalign{\medskip}
{\frac{\delta \beta_2 \eta I_T}{N}}
&{\frac{\delta \beta_2 I_T}{N}}
&{\frac{\delta \beta_2 I_T}{N}}
&{\frac{\delta \beta_2 I_T}{N}}
&{\frac{\delta \beta_2 \eta  I_T}{N}}\\
\noalign{\medskip}
{\frac{\beta_2 \eta R_T}{N}}
&{\frac{\beta_2 R_T}{N}}
& \left({\frac{\beta^{'}_1 \beta_1}{N}}
+{\frac{\beta_2}{N}} \right) R_T
&{\frac{\beta_2 R_T}{N}}
&\left(\frac{\beta_1^{'} \beta_1}{N}+ {\frac{\beta_2 \eta }{N}} \right)R_T\\
\noalign{\medskip}
-\alpha_1&0&
{\frac{\psi \beta_1 I_H}{N}}&0&{\frac{\psi \beta_1 I_H}{N}}\\
\noalign{\medskip}
\alpha_1+\mu+d_A&0&0&0&0\\
\noalign{\medskip}
0&k_2+\tau_4 +\mu&0&0&0\\
\noalign{\medskip}
0&-k_2&\rho_2+\tau_3+\mu+d_T&0&-\alpha_2\\
\noalign{\medskip}
0&-\tau_4&-\tau_3+{\frac{\beta^{'}_2 \beta_1 R_{TH}}{N}}&\beta^{'}_2
\lambda_T+\rho_3+\mu&\frac{\beta_2^{'} \beta_1 R_{TH}}{N}\\
\noalign{\medskip}
0&0&-\rho_2&-\rho_3&\alpha_2+d_{TA}+\mu\end {array} \right].
\end{equation*}
}}
The dominant eigenvalues of the matrix $F \cdot V^{-1}$ are
\begin{equation*}
\begin{split}
R_1&=\frac{\Lambda}{N\mu}\left(\frac{\beta_1}{d_T+\mu+\tau_2} \right)
\left(\frac{k_1}{k_1+\tau_1+\mu } \right) \, ,\\
R_2&=\frac{\Lambda}{N\mu} \beta_2 \left(\frac{\mu+\alpha_1+d_A
+\eta\,\rho_1 }{\mu\,\alpha_1+ (\mu + \rho_1)(\mu + d_A)}\right)\, .
\end{split}
\end{equation*}
Thus, the basic reproduction number $R_0$
of the model \eqref{model:TB:HIV} is given by
\begin{equation*}
R_0 = \max \{ R_1, R_2 \} \, .
\end{equation*}
Note that $R_1$ is the basic reproduction number of the model \eqref{model:TB:HIV}
with $I_T = A = L_{TH} = I_{TH} = R_{TH} = A_T = 0$ (only TB model), and $R_2$
is the basic reproduction number of the model \eqref{model:TB:HIV} with
$L_T = I_T = R_T = L_{TH} = I_{TH} = R_{TH} = A_T = 0$ (only HIV-AIDS model).


\section{Appendix: Proof of Theorem~\ref{theo:stab:equil}}
\label{A.2}

In this Appendix we provide details of the proof of Theorem~\ref{theo:stab:equil}.

\paragraph{Local asymptotical stability of the disease-free equilibrium $\Sigma_0$.}

Following Theorem~2 of \cite{van:den:Driessche:2002}, the disease-free equilibrium,
$\Sigma_0$, is locally asymptotically stable if all the eigenvalues of the Jacobian
matrix of the system \eqref{model:TB:HIV}, here denoted by $M_T\left(\Sigma_0\right)$,
computed at the disease free equilibrium $\Sigma_0$, given by \eqref{eq:DFE:model:TBHIV},
have negative real parts.

The Jacobian matrix of the system \eqref{model:TB:HIV}
at disease free equilibrium $\Sigma_0$ is given by
\begin{equation*}
M_T\left(\Sigma_0\right) = \left[
\begin{array}{cc}
M_{T1}\left(\Sigma_0\right) \quad M_{T_2}\left(\Sigma_0\right)
\end{array}
\right]
\end{equation*}
with
\begin{equation*}
M_{T1}\left(\Sigma_0\right)=
{\footnotesize{
\left[ \begin {array}{ccccc}
-\mu&0&-{\frac{\beta_1 \Lambda}{\mu N}}&0&
-{\frac{\beta_2 \Lambda}{\mu N}}\\ \noalign{\medskip}
0&-d_1&{\frac{\beta_1 \Lambda}{\mu N}}&0&0\\ \noalign{\medskip}
0&k_1&-d_2&0&0\\ \noalign{\medskip}
0&\tau_1&\tau_2&-\mu&0\\ \noalign{\medskip}
0&0&0&0&{\frac{\beta_2 \Lambda}{\mu N}}-d_3\\ \noalign{\medskip}
0&0&0&0&\rho_1\\ \noalign{\medskip}
0&0&0&0&0\\ \noalign{\medskip}
0&0&0&0&0\\ \noalign{\medskip}
0&0&0&0&0\\ \noalign{\medskip}
0&0&0&0&0
\end {array} \right]
}}
\end{equation*}
and
\begin{equation*}
M_{T2}\left(\Sigma_0\right)=
{\tiny{
\left[ \begin {array}{cccccccccc}
-{\frac{\beta_2 \eta \Lambda}{\mu N}}&
-{\frac{\beta_2 \Lambda}{\mu N}}&
-{\frac {\beta_1 \Lambda}{\mu N}}-{\frac{\beta_2 \Lambda}{\mu N}}&
-{\frac {\beta_2 \Lambda}{\mu N}}&
-{\frac{\beta_1 \Lambda}{\mu N}}-{\frac{\beta_2 \eta  \Lambda}{\mu N}}\\
\noalign{\medskip}
0&0
&{\frac{\beta_1 \Lambda}{\mu N}}&0&{\frac{\beta_1 \Lambda}{\mu N}}\\
\noalign{\medskip}
0&0&0&0&0\\
\noalign{\medskip}
0&0&0&0&0\\
\noalign{\medskip}
{\frac{\beta_2 \eta \Lambda}{\mu N}}+\alpha_1
&{\frac {\beta_2 \Lambda}{\mu N}}
&{\frac{\beta_2 \Lambda}{\mu N}}
&{\frac{\beta_2 \Lambda}{\mu N}}
&{\frac {\beta_2 \eta  \Lambda}{\mu N}}\\
\noalign{\medskip}
-d_4&0&0&0&0\\ \noalign{\medskip}
0&-d_5&0&0&0\\ \noalign{\medskip}
0&k_2&-d_6&0&\alpha_2\\ \noalign{\medskip}
0&\tau_4&\tau_3&-d_7&0\\ \noalign{\medskip}
0&0&\rho_2&\rho_3&-d_8
\end {array} \right] ,
}}
\end{equation*}
where
$d_1 = k_1 +\tau_1+\mu$; $d_2 = \tau_2+\mu+d_T$; $d_3 = \rho_1+\mu$;
$d_4 = \alpha_1+\mu+d_A$; $d_5 = k_2+\mu + \tau_4$;
$d_6 = \rho_2+\tau_3+\mu+d_T$; $d_7 = \rho_3+\mu$;
$d_8 = \alpha_2+ d_{TA}+\mu$. One has
\begin{equation*}
trace\left[ M_T\left(\Sigma_0\right) \right] =
-2 \mu - (d_1 + d_2 + d_3 + d_4 + d_5 + d_6 + d_7 + d_8) < 0
\end{equation*}
and
\begin{equation*}
\begin{split}
det \left[ M_T\left(\Sigma_0\right) \right] &
= \frac{1}{N^2} (
d_5\, ( d_6 \, d_7 + d_T(\alpha_2 + \mu)\,d_7 + \alpha_2 \mu \,d_6 + d_T d_{TA}\, d_7 )\\
&\times (N \mu (\alpha_1 \mu + (\mu + \rho_1)(\mu + d_A))
-\beta_2 \Lambda(\alpha_1 + \mu + d_A + \rho_1 \eta))\\
&\times (N \mu (d_T + \mu + \tau_2)
(k_1 + \tau_1 + \mu)-k_1 \beta_1 \Lambda) > 0
\end{split}
\end{equation*}
for
\begin{equation*}
R_1 = \frac{\Lambda}{N\mu}\left(\frac{\beta_1}{d_T+\mu+\tau_2} \right)
\left(\frac{k_1}{k_1+\tau_1+\mu } \right) < 1
\end{equation*}
and
\begin{equation*}
R_2 =\frac{\Lambda}{N\mu} \beta_2 \left(\frac{\mu+\alpha_1+d_A
+\eta\,\rho_1 }{\mu\,\alpha_1+ (\mu + \rho_1)(\mu + d_A)}\right) < 1 \, .
\end{equation*}
We have just proved that the disease free equilibrium $\Sigma_0$ of model
\eqref{model:TB:HIV} is locally asymptotically stable if $R_0 < 1$,
and unstable if either $R_i > 1$, $i=1, 2$.


\bigskip


\paragraph{Global asymptotical stability of the disease-free equilibrium $\Sigma_0$.}

For convenience, let us rewrite the model system \eqref{model:TB:HIV} as
\begin{equation}
\label{mod:append}
\begin{split}
&\frac{dX}{dt} = F(X, Z) \, , \\
&\frac{dZ}{dt} = G(X, Z)\, , \quad G(X, 0) = 0 \, ,
\end{split}
\end{equation}
where $X = (S, R_T)$ and $Z = (L_T, I_T, I_H, A, L_{TH}, I_{TH}, R_{TH}, A_T)$,
with $X \in \R^2_+$ denoting (its components) the number of uninfected individuals
and $Z \in \R^8_+$ denoting (its components) the number of infected
individuals including the latent and infectious.

The disease-free equilibrium is denoted by
\begin{equation*}
E_0 = (X_0, 0)\, , \quad \text{where} \, \, X_0
= \left(\frac{\Lambda}{\mu}, 0  \right) \, .
\end{equation*}
Following \cite{Bhunu:BMB:2009:HIV:TB}, if
\begin{itemize}
\item[(H1)] $E_0$ is globally asymptotically stable for $\frac{dX}{dt} = F(X, 0)$,
\item[(H2)] $\hat{G}(X, Z) \geq 0$ for $(X, Z) \in \Omega$, where $G(X, Z)
= AZ - \hat{G}(X, Z)$, $A = D_Z G(E_0, 0)$ is a Metzler matrix and
$\Omega$ is given by \eqref{eq:feasible:region},
\end{itemize}
then the fixed point $E_0 = (X_0, 0)$ is a globally asymptotically
stable equilibrium of system \eqref{mod:append}. We have
\begin{equation*}
\frac{dX}{dt} = F(X, Z)
= \left[
\begin{array}{c}
\Lambda - \lambda_T S - \lambda_H S - \mu S \\[0.2 cm]
\tau_1 L_T + \tau_2 I_T - (\beta^{'}_1 \lambda_T + \lambda_H + \mu) R_T
\end{array} \right] \, ,
\end{equation*}
\begin{equation*}
F(X, 0) = \left[ \begin{array}{c}
\Lambda - \mu S \\[0.2 cm]
-\mu R_T
\end{array} \right] \, ,
\end{equation*}
\begin{equation*}
\frac{dZ}{dt} = G(X, Z) = \left[
\begin{array}{c}
\lambda_T S+ \beta^{'}_1 \lambda_T R_T- (k_1 + \tau_1 + \mu)L_T,\\[0.2 cm]
k_1 L_T - (\tau_2 +d_T +\mu + \delta \lambda_H)I_T, \\[0.2 cm]
\lambda_H S - (\rho_1 + \psi \lambda_T + \mu)I_H + \alpha_1 A + \lambda_H R_T \\[0.2 cm]
\rho_1 I_H - \alpha_1 A - (\mu + d_A) A\\[0.2 cm]
\beta^{'}_2 \lambda_T R_{TH} - (k_2 + \tau_4 + \mu) L_{TH}\\[0.2 cm]
\delta \lambda_H I_T + \psi \lambda_T I_H + \alpha_2 A_T+ k_2 L_{TH}
- (\tau_3 + \rho_2 + \mu + d_T)I_{TH}\\[0.2 cm]
\tau_3 I_{TH} + \tau_4 L_{TH} - (\beta^{'}_2 \lambda_T  + \rho_3 + \mu)R_{TH} \\[0.2 cm]
\rho_2 I_{TH} + \rho_3 R_{TH} -(\alpha_2 + \mu + d_{TA})A_T \, ,
\end{array}
\right],
\end{equation*}
and $G(X, 0) = 0$. Thus,
\begin{equation*}
\frac{dX}{dt} = F(X, 0) = \left[
\begin{array}{c}
\Lambda - \mu S \\[0.2 cm]
- \mu R_T
\end{array} \right] \, ,
\end{equation*}
\begin{equation*}
A = D_Z G(X_0, 0) = \left[ \begin {array}{cc} D_1 & D_2
\end{array} \right]
\end{equation*}
with
\begin{equation*}
D_1 = {\footnotesize{  \left[ \begin {array}{cccc}
-k_1-\tau_1-\mu & \frac{\beta_1 \Lambda}{\mu N}&0&0
\\ \noalign{\medskip}
k_1&-\tau_2-\mu- d_T&0&0\\ \noalign{\medskip}
0&0&\frac {\beta_2 \Lambda}{\mu N}-\rho_1-\mu&
\frac{\beta_2 \eta \Lambda}{\mu N}+\alpha_1\\ \noalign{\medskip}
0&0&\rho_1&-\alpha_1-\mu- d_A\\ \noalign{\medskip}
0&0&0&0\\ \noalign{\medskip}
0&0&0&0\\ \noalign{\medskip}
0&0&0&0\\ \noalign{\medskip}
0&0&0&0
\end {array} \right] }} \, ,
\end{equation*}
\begin{equation*}
D_2 = {\footnotesize{  \left[ \begin {array}{cccccccc}
0&{\frac{\beta_1 \Lambda}{\mu N}}&0&{\frac{\beta_1 \Lambda}{\mu N}}\\
\noalign{\medskip}
0&0&0&0\\ \noalign{\medskip}
{\frac{\beta_2 \Lambda}{\mu N}}&
{\frac{\beta_2 \Lambda}{\mu N}}&{\frac{\beta_2 \Lambda}{\mu N}}
&{\frac{\beta_2 \eta  \Lambda}{\mu N}}\\ \noalign{\medskip}
0&0&0&0\\ \noalign{\medskip}
-k_2-\tau_4 -\mu &0&0&0\\ \noalign{\medskip}
k_2&-\rho_2-\tau_3-\mu - d_T&0&\alpha_2\\ \noalign{\medskip}
\tau_4&\tau_3&-\rho_3-\mu&0\\ \noalign{\medskip}
0&\rho_2&\rho_3&-\alpha_2-d_{TA}-\mu
\end {array} \right] }}
\end{equation*}
and
\begin{equation}
\label{eq:hat:G}
\hat{G}(X, Z)
= {\footnotesize{  \left[ \begin {array}{c}
\lambda_T \left(\frac{\Lambda}{\mu} - S - \beta_1^{'} R_T  \right)\\
\noalign{\medskip}
-\delta \lambda_H I_T\\
\noalign{\medskip}
\lambda_H \left( \frac{\Lambda}{\mu} - S - R_T - \psi I_H \right)\\
\noalign{\medskip}
0\\ \noalign{\medskip}
-\beta_2^{'} \lambda_T R_{TH}\\ \noalign{\medskip}
-\left(\delta \lambda_H I_T + \psi \lambda_T I_H \right)\\
\noalign{\medskip}
\beta_2^{'} \lambda_T R_{TH}\\
\noalign{\medskip}
0
\end {array} \right] }} \, .
\end{equation}
From \eqref{eq:hat:G} the condition (H2) is not satisfied, since $\hat{G}(X, Z) \geq 0$
is not true. Therefore, the disease-free equilibrium $E_0$ may not be globally
asymptotically stable. Following \cite{CChavez_TB_Exog_TPB_2000}, the backward
bifurcation occurs at $R_0 = 1$ and the double endemic equilibria can be supported
for $R_c < R_0 < 1$, where $R_c$ is a positive constant.


\bigskip


\paragraph{Existence and stability of HIV-AIDS free equilibrium $\Sigma_T$.}

The expressions for $S^\diamond$, $L_T^\diamond$, $I_T^\diamond$ and $R_T^\diamond$
are obtained if we consider a sub-model of \eqref{model:TB:HIV} for which
$I_H = A = L_{TH} = I_{TH} = R_{TH} = A_T = 0$ and the total population $N$
is given by $N_T = S + L_T + I_T + R_T$. The basic reproduction number of this
submodel is given by $R_1$ \eqref{eq:R1}. The existence, uniqueness and local
asymptotic stability of $\Sigma_T$ is proven in \cite[Theorem~1]{Castillo_Chavez_1997}.


\bigskip


\paragraph{Existence and stability of TB free equilibrium $\Sigma_H$.}

To prove the existence of $\Sigma_T$, consider the sub-model of \eqref{model:TB:HIV}
for which $L_T = I_T = R_T = L_{TH} = I_{TH} = R_{TH} = A_T = 0$
and the total population $N_H$ is given by $N_H = S + I_H + A$.
The equations of this submodel are
\begin{equation}
\label{submodel:HIV}
\begin{cases}
\dot{S}(t) = \Lambda - \lambda_H S(t) - \mu S(t)\\[0.2 cm]
\dot{I}_H(t) = \lambda_H S(t) - (\rho_1 + \mu)I_H(t) + \alpha_1 A(t)\\[0.2 cm]
\dot{A}(t) = \rho_1 I_H(t) - \alpha_1 A(t) - (\mu + d_A)A \, ,
\end{cases}
\end{equation}
where $\lambda_H = \beta_2 \frac{I_H + \eta A}{N_H}$.
Setting the right hand sides of submodel \eqref{submodel:HIV} to zero,
we obtain the endemic equilibrium $\Sigma_H^\star = (S^*, I_H^*, A^*)$ given by
\begin{equation*}
S^\star = \frac{\Lambda}{\mu R_2}\, , \quad I_H^\star
= (R_2 - 1)\frac{\mu N_H (\alpha_1 + d_A + \mu)}{\beta_2 (\alpha_1
+ d_A + \mu + \eta \rho_1)}\, , \, \, \quad  A^\star
= (R_2 - 1) \frac{\rho_1 \mu N_H}{\beta_2 (\alpha_1 + d_A + \mu + \eta \rho_1)},
\end{equation*}
where $I_H^\star > 0$ and $A^\star > 0$, whenever $R_2 > 1$.

In what follows we prove the local asymptotic stability of the endemic equilibrium
$\Sigma_H^\star$, using the center manifold theory \cite{Carr:1981},
as described in \cite[Theorem~4.1]{CChavez_Song_2004}
(see also \cite{van:den:Driessche:2002}), considering ART treatment.
The basic reproduction number of this sub-model $R_2$ is given by \eqref{eq:R2}.
Chose as bifurcation parameter, $\beta^*$, by solving for $\beta_2$ from $R_2 = 1$:
\begin{equation*}
\beta^* = {\frac {\mu \alpha_1 + (\mu + \rho_1)(\mu + d_A)}{
\alpha+ d_A +\mu+\eta \rho}} \, .
\end{equation*}
The submodel \eqref{submodel:HIV} has a disease free equilibrium given by
$\Sigma_{H0}^*=(x_{10}, x_{20}, x_{30})= \left(\frac{\Lambda}{\mu}, 0, 0 \right)$.

The Jacobian of the system \eqref{submodel:HIV}, evaluated at
$\Sigma_{H0}^*$ and with $\beta_2 = \beta^*$, is given by
\begin{equation}
\label{Jacob:betaast}
J(\Sigma_{H0}^*) = \left[ \begin {array}{ccc}
-\mu&- \beta_2&-\beta_2 \eta\\
\noalign{\medskip}
0&\beta_2-\rho-\mu & \beta_2 \eta + \alpha\\
\noalign{\medskip}
0&\rho&-\alpha- d_A - \mu
\end {array} \right] \, .
\end{equation}
The eigenvalues of the linearized system \eqref{Jacob:betaast} are
$$
\lambda_1 = 0, \quad \lambda_2 = -\mu \quad \text{and} \quad \lambda_3
= -\frac{ \eta \rho (2\mu^2 + \rho + d_A + \alpha) + d_A(2 \alpha + 2 \mu +d_A)
+\rho \alpha+ (\mu + \alpha)^2}{\alpha+ d_A+\mu+ \eta \rho} \, .
$$
We observe that there is a simple eigenvalue with zero real part and the other
two eigenvalues have negative real part. Thus, the system \eqref{submodel:HIV},
with $\beta_2 = \beta^*$, has a hyperbolic equilibrium point and the center
manifold theory \cite{Carr:1981} can be used to analyze the dynamics
of the submodel \eqref{submodel:HIV} near $\beta_2 = \beta^*$.

The Jacobian $J(\Sigma_{H0}^*)$ at $\beta_2 = \beta^*$ has a right eigenvector
(associated with the zero eigenvalue) given by $w = [w_1, w_2, w_3]^T$, where
\begin{equation*}
\begin{split}
w_1&=-\frac{\left( \mu \alpha_1+ (\mu + \rho_1)(\mu + d_A) \right) w_3}{\rho_1 \mu}\, ,\\
w_2&=\frac { \left( \alpha_1 + d_A+\mu \right) w_3}{\rho_1}\, ,\\
w_3&= w_3 > 0 \, .
\end{split}
\end{equation*}
Further, $J(\Sigma_{H0}^*)$ for $\beta_2 = \beta^*$ has a left eigenvector
$v = [v_1, v_2, v_3]$ (associated with the zero eigenvalue), where
\begin{equation*}
\begin{split}
v_1&=0\, ,\\
v_2&=\frac{v_3 \left( \alpha_1+ d_A +\mu+\eta \rho_1 \right)}{
\alpha_1 +\eta \rho_1 +\mu \eta}\, ,\\
v_3&= v_3 > 0 \, .
\end{split}
\end{equation*}
To apply Theorem~4.1 in \cite{CChavez_Song_2004} it is convenient to let
$f_k$ represent the right-hand side of the $k$th equation of the system
\eqref{submodel:HIV} and let $x_k$ be the state variable whose derivative
is given by the $k$th equation for $k = 1, 2, 3$. The local stability
near the bifurcation point $\beta_2 = \beta^*$ is then determined by the
signs of two associated constants, denoted by $a$ and $b$, defined (respectively) by
\begin{equation*}
a = \sum_{k, i, j=1}^3 \, v_k w_i w_j \frac{\partial^2 f_k}{\partial x_i \partial x_j}(0, 0)
\quad \text{and} \quad b = \sum_{k, i =1}^3 \, v_k w_i
\frac{\partial^2 f_k}{\partial x_i \partial \phi}(0, 0)
\end{equation*}
with $\phi = \beta_2 - \beta^*$.

For the system \eqref{submodel:HIV}, the associated partial
derivatives at the disease free equilibrium $\Sigma_{H0}$ are given by
\begin{equation*}
\begin{split}
\frac{\partial^2 f_1}{\partial x_2^2} &= \frac{2\beta^*\mu}{\Lambda}\, ,  \quad
\frac{\partial^2 f_1}{\partial x_2 \partial x_3} = \frac{\beta^*\mu(1+\eta)}{\Lambda}\, , \quad
\frac{\partial^2 f_1}{\partial x_3^2} = \frac{2\beta^*\mu \eta}{\Lambda} \, ,\\
\frac{\partial^2 f_2}{\partial x_2^2} &= \frac{-2\beta^*\mu}{\Lambda}\, ,  \quad
\frac{\partial^2 f_2}{\partial x_2 \partial x_3} = \frac{-\beta^*\mu(1+\eta)}{\Lambda}\, , \quad
\frac{\partial^2 f_2}{\partial x_3^2} = \frac{-2\beta^*\mu \eta}{\Lambda} \, .
\end{split}
\end{equation*}
It follows from the above expressions that
\begin{equation*}
a=-\frac {v_3 w_3^2 \beta^* \mu \left(k_1 +\mu+\eta \rho_1 \right)
\left( 2 k_1^2  + 4 \mu k_1 + 2 \mu^2 + \rho_1 (\alpha_1
+ \eta (\alpha_1 + \mu + 2 \rho_1) + d_A(1 + \eta) + \mu)
\right)}{\rho_1^2\Lambda \left( \alpha_1
+\eta \rho_1 +\mu \eta \right)} < 0
\end{equation*}
with $k_1 = \alpha_1+ d_A $.

For the sign of $b$, it can be shown that
the associated non-vanishing partial derivatives are
\begin{equation*}
\begin{split}
\frac{\partial^2 f_1}{\partial x_2 \partial \beta^*} &= -1\, ,  \quad
\frac{\partial^2 f_1}{\partial x_3 \partial \beta^*} = -\eta \, , \\
\frac{\partial^2 f_2}{\partial x_2 \partial \beta^*} &= 1 \, , \quad
\frac{\partial^2 f_2}{\partial x_3 \partial \beta^*} = \eta \, .
\end{split}
\end{equation*}
It also follows from the above expressions that
\begin{equation*}
b=\frac{ v_3 w_3 \left( k_1 +\mu+\eta \rho_1 \right)
 \left( k_1 +\mu \right) }{ \left( \alpha_1 +\eta
\rho_1 +\mu \eta \right) \rho_1}+ \frac {\eta v_3 w_3
\left( k_1 +\mu+\eta \rho_1 \right)}{\alpha_1 +\eta \rho_1 + \mu \eta} > 0 \, .
\end{equation*}
Thus, $a< 0$ and $b > 0$. Using Theorem~4.1 of \cite{CChavez_Song_2004},
the endemic equilibrium $\Sigma_H^\star$
is locally asymptotically stable for $R_2$ near 1.


\section*{Acknowledgments}

This work was supported by Portuguese funds through the
\emph{Center for Research and Development in Mathematics and Applications}
(CIDMA, University of Aveiro), and \emph{The Portuguese Foundation for Science and Technology}
(``FCT --- Funda\c{c}\~{a}o para a Ci\^{e}ncia e a Tecnologia''),
within project PEst-OE/MAT/UI4106/2014.
Silva was also supported by FCT through
the post-doc fellowship SFRH/BPD/72061/2010;
Torres by the FCT project PTDC/EEI-AUT/1450/2012,
co-financed by FEDER under POFC-QREN with COMPETE
reference FCOMP-01-0124-FEDER-028894.




\begin{thebibliography}{99}

\bibitem{Bhunu:BMB:2009:HIV:TB}
C.~P.~Bhunu, W.~Garira and Z.~Mukandavire,
\textit{Modeling HIV/AIDS and tuberculosis coinfection},
Bul. Math. Biol. 71, 1745--1780 (2009).

\bibitem{Carr:1981}
J.~Carr,
\textit{Applications centre manifold theory},
Springer-Verlag, New-York (1981).

\bibitem{Castillo_Chavez_1997}
C.~Castillo-Chavez and Z.~Feng,
\textit{To treat or not to treat: The case of tuberculosis},
J. Math. Biol. 35, no.~6, 629--656 (1997).

\bibitem{CChavez_Song_2004}
C.~Castillo-Chavez and B.~Song,
\textit{Dynamical models of tuberculosis and their applications},
Math. Biosc. Engrg. 1, no.~2, 361--404 (2004).

\bibitem{art:viral:load}
P.~W.~David, G.~L.~Matthew, E.~G.~Andrew, A.~C.~David and M.~K.~John,
\textit{Relation between HIV viral load and infectiousness:
A model-based analysis},
The Lancet 372, no.~9635, 314--320 (2008).

\bibitem{CChavez_TB_Exog_TPB_2000}
Z.~Feng, C.~Castillo-Chavez and A.~F.~Capurro,
\textit{A model for tuberculosis with exogenous reinfection},
Theor. Pop. Biology 57, 235--247 (2000).

\bibitem{Getahun:etall:CID2010}
H.~Getahun, C.~Gunneberg, R.~Granich and P.~ Nunn,
\textit{HIV infection-associated tuberculosis:
The epidemiology and the response},
Clin. Infect. Dis. 50 (Suppl 3), S201--S207 (2010).

\bibitem{Hyman:MathBio:1999}
J.~M.~Hyman, J.~Li and E.~A.~Stanley,
\textit{The differential infectivity
and staged progression models for the transmission of HIV},
Math. Biosci. 155, 77--109 (1999).

\bibitem{Kirschner:TB:HIV:1999}
D.~Kirschner,
\textit{Dynamics of co-infection with M. tuberculosis and HIV-1},
Theor. Pop. Biol. 55, 94--109 (1999).

\bibitem{Kwan_Ernst_HIV_TB_Syndemic}
C.~K.~Kwan and J.~D.~Ernst,
\textit{HIV and tuberculosis: A deadly human syndemic},
Clin. Microbiol. Rev. 24, no.~2, 351--376 (2011).

\bibitem{Naresh:TB:HIV:2005}
R.~Naresh and A.~Tripathi,
\textit{Modelling and analysis of HIV-TB
co-infection in a variable size population},
Math. Model. Anal. 10, 275--286 (2005).

\bibitem{CChavez_TB_HIV_2009}
L.~W.~Roeger, Z.~Feng and C.~Castillo-Chavez,
\textit{Modeling TB and HIV co-infections},
Math. Biosc. and Eng. 6, no.~4, 815--837 (2009).

\bibitem{Song:TB:HIV:2008}
O.~Sharomi, C.N.~Podder, A.B.~Gumel and B.~Song,
\textit{Mathematical analysis of the transmission dynamics
of HIV/TB coinfection in the presence of treatment},
Math. Biosc. Eng. 5, no.~1, 145--174 (2008).

\bibitem{UNAIDS_report_2013}
UNAIDS,
\textit{Global report: UNAIDS report on the global AIDS epidemic 2013},
Geneva, World Health Organization (2013).

\bibitem{van:den:Driessche:2002}
P.~van den Driessche and J.~Watmough,
\textit{Reproduction numbers and subthreshold endemic equilibria
for compartmental models of disease transmission},
Math. Biosc. 180, 29--48 (2002).

\bibitem{TB_WHO_report_2013}
WHO,
\textit{Global tuberculosis report 2013},
Geneva, World Health Organization (2013).

\bibitem{USAID:TB:HIV}
\url{http://www.usaid.gov/news-information/fact-sheets/twin-epidemics-hiv-and-tb-co-infection}

\bibitem{wiki:HIV:progression}
\url{http://en.wikipedia.org/wiki/HIV_disease_progression_rates}

\bibitem{site:HIV:progression:rate}
\url{http://hivinsite.ucsf.edu/InSite?page=kb-03-01-04}

\end{thebibliography}
\end{document}